\documentclass{amsart}
\numberwithin{equation}{section}

\usepackage{amsmath,amssymb,mathtools}
\usepackage{booktabs,longtable,array}
\usepackage{enumitem}
\usepackage{microtype}

\newcommand{\CD}[1]{{}^{\mathrm C}\!D_t^{#1}}
\newcommand{\RL}[1]{{}^{\mathrm{RL}}\!D_t^{#1}}
\newcommand{\Ii}[1]{I_t^{#1}}
\begin{document}

\author{Ravshan Ashurov}
\address{Ashurov R: Institute of Mathematics of Uzbekistan,
Tashkent, Uzbekistan}
\email{ashurovr@gmail.com}

\author{Nigora Qutlimurotova}
\address{Qutlimurotova N:National University of Uzbekistan,  Tashkent, Uzbekistan}
\email{n.kutlimurotova@nuu.uz}

\author{Ilyoskhuja Sulaymonov}
\address{Sulaymonov I: National University of Uzbekistan,  Tashkent, Uzbekistan}
\email{sulaymonov\_i@nuu.uz}

\small

\title[Fractional Differential Equations in Chemical Modeling ...] {Fractional Differential Equations in Chemical Modeling: From Memory and Anomalous Transport to Mathematically Justified Models}

\begin{abstract}

Fractional differential equations are increasingly used to model memory and anomalous transport in chemistry and related fields. However, the formal replacement of an integer-order derivative by a fractional one does not by itself justify fractional dynamics. The essential issue is whether nonlocality can be derived from a physical, stochastic, geometric, or reduction mechanism and whether it leads to experimentally testable consequences.

This review, written primarily for mathematicians and specialists in fractional analysis, discusses representative applications in electrochemistry, adsorption and sorption, reaction--subdiffusion, chemical and enzyme kinetics, polymer relaxation, and atmospheric chemistry. Particular attention is paid to distinguishing mechanism-based fractional models from phenomenological fractionalizations.

Mathematical questions of well-posedness, nonnegativity, conservation laws, dimensional consistency, parameter identifiability, and numerical convergence are emphasized. Experimental electrochemical examples illustrate the need to distinguish genuine fractional effects from model mismatch or non-identifiability. Several open problems connecting fractional analysis with experimentally grounded chemical modeling are formulated.

\end{abstract}

\maketitle

\section{Introduction}

Differential equations are among the principal mathematical tools of chemistry. Reaction kinetics, diffusion of dissolved substances, transport to electrodes, adsorption and desorption at interfaces, dielectric relaxation, and atmospheric chemical processes are traditionally described by ordinary or partial differential equations of integer order.

The classical approach is highly successful in many situations. For example, an irreversible first-order chemical transformation can be represented by
\begin{equation}\label{eq:first}
\frac{dC}{dt}=-kC,
\end{equation}
where $C(t)$ is the concentration of the reacting species and $k>0$ is the reaction-rate constant. The term \emph{first-order reaction} means that the disappearance rate is proportional to the current concentration. Thus, for constant $k$, $C(t)=C(0)e^{-kt}$, and the relaxation is governed by one exponential and one characteristic time $\tau=1/k$.

Similarly, normal diffusion in a homogeneous medium is described by Fick's law together with the local conservation law,
\begin{equation}\label{eq:fick}
J=-D\nabla C,\qquad \frac{\partial C}{\partial t}=D\Delta C,
\end{equation}
where $J$ is the molecular flux and $D$ is the diffusion coefficient. Here $\nabla C$ is the spatial concentration gradient, while $\Delta=\nabla\cdot\nabla$ is the Laplace operator; in Cartesian coordinates $\Delta C=\sum_j\partial^2C/\partial x_j^2$.

Model \eqref{eq:first} contains a single exponential relaxation mode. Many chemical and physicochemical systems, however, exhibit broad relaxation spectra, trapping, barriers, pore-scale heterogeneity, interfacial effects, or molecular configurations that reorganize over many time scales. Likewise, transport may violate the classical scaling of normal diffusion. In such cases the current rate of change may depend not only on the present state but also on the previous history of the process.

The main purpose of this article is to identify chemical problems in which nonlocal mathematical models arise for a concrete reason and in which fractional analysis produces consequences that cannot be reduced to a cosmetic change of notation. The aim is not to catalogue every paper containing a fractional derivative, but to separate several conceptually different mechanisms by which a fractional operator can enter a chemical model.

\begin{figure}[t]
\centering
\includegraphics[width=0.95\textwidth]{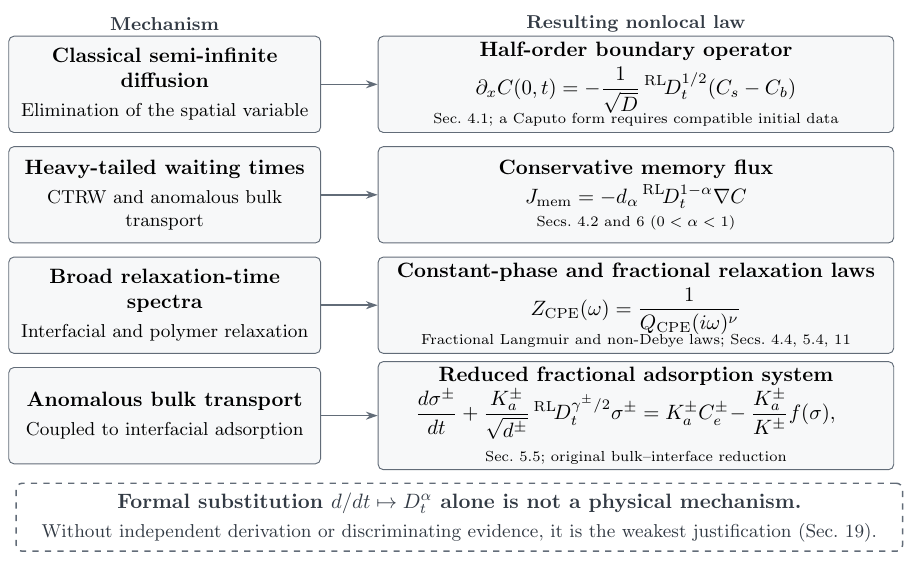}
\caption{Mechanisms that generate nonlocal operators in the chemical models discussed in this article, and the operators they lead to. The dashed box marks formal fractionalization without an identified mechanism.}
\label{fig:scheme}
\end{figure}

This distinction is fundamental. If one starts from $du/dt=F(u)$ and formally replaces it by
\[
D_t^\alpha u=F(u),\qquad 0<\alpha<1,
\]
this still does not explain why the chemical system should possess fractional dynamics. A convincing model should answer the more fundamental question
\begin{center}
\emph{What physical, stochastic, geometric mechanism, or reduction procedure generates the nonlocality?}
\end{center}
Possible answers include heavy-tailed waiting times, anomalous diffusion, a distribution of relaxation times, unresolved spatial heterogeneity, elimination of hidden variables, or a constitutive law with memory. The modeling chain is therefore
\begin{align*}
&\text{chemical phenomenon}\to\text{classical model}\to\text{observable limitation}\\
&\quad\to\text{memory or anomalous-transport mechanism}\to\text{nonlocal model}\to\text{fractional operator}.
\end{align*}
This logic is summarised in Fig.~\ref{fig:scheme}.

A second purpose is methodological. Fractional versions of Brusselator systems \cite{asiedu2025}, pollutant transport in water \cite{patel2023}, soil and wetland carbon dynamics \cite{bohaienko2025} and even quantum chemistry \cite{mostafanejad2021} have been proposed. A substantial part of this literature concentrates on constructing approximate solutions by series, transform--decomposition, operational-matrix or wavelet techniques \cite{shiri2025}, and contains numerous numerical computations. Numerical computation is indispensable, but a smooth numerical curve is not a substitute for a well-posed continuous problem. If uniqueness has not been established, a numerical algorithm may converge to one among several admissible solutions. Even when the continuous problem is uniquely solvable, convergence of the discrete approximation to that solution is a separate mathematical question.

\textit{Scope and literature-selection criteria.} 
The literature was selected purposively rather than by a systematic-review protocol. A work was included in the detailed discussion when it satisfied at least one of three criteria: (i) the nonlocal operator was derived from a transport, stochastic, constitutive, or reduction mechanism; (ii) a fractional model was calibrated against or compared with real chemical or physicochemical measurements; or (iii) the work contained a mathematically substantial result---for example, a well-posedness theorem, qualitative analysis, or a convergent numerical method---for a chemically motivated fractional model. Purely formal replacements of an integer derivative by a fractional one are mentioned only when they represent a broader trend or when they are useful for contrasting levels of evidence. These three categories---mechanism-derived, experimentally confronted, and mathematically or computationally motivated---are kept distinct in the discussion and are identified explicitly in Table~\ref{tab:status}.

\section{Chemical concepts needed for the mathematical models}

\textit{ 1. Concentration, reaction rate, and reaction order.}
The concentration of a chemical species measures the amount of that species per unit volume. We write it as $C(t)$ in a spatially homogeneous model and as $C(x,t)$ when spatial transport is present.

For a simple transformation $A\longrightarrow B$, the reaction rate quantifies the rate at which $A$ disappears or $B$ is formed. In the first-order law $-dC_A/dt=kC_A$ the rate is linear in the reactant concentration. More complicated reaction mechanisms involve several species and nonlinear mass-action terms, so chemical kinetics naturally leads to systems of nonlinear differential equations.

The word \emph{order} in chemical kinetics refers to the exponent with which a concentration enters the empirical or mechanistic rate law. It should not be confused with the order of a fractional derivative. These are mathematically different notions and may coexist in the same model.

\textit{2. Diffusion, flux, and mean square displacement.}
Diffusion is the net molecular transport generated by spatial nonuniformity of concentration. The flux $J(x,t)$ is the amount of substance crossing a unit area per unit time. In the classical Fickian model, $J=-D\nabla C$, and combining this constitutive law with local mass conservation $\partial_tC=-\nabla\cdot J$ gives \eqref{eq:fick}.

For normal diffusion, the mean square displacement satisfies $\langle|X(t)-X(0)|^2\rangle\sim t$. The mean square displacement is the ensemble average of the squared distance travelled by a particle from its initial position and measures the characteristic spatial scale of spreading. If instead
\begin{equation}\label{eq:msd}
\langle|X(t)-X(0)|^2\rangle\sim t^\gamma,\qquad\gamma\neq1,
\end{equation}
then the transport is called anomalous. Values $0<\gamma<1$ correspond to subdiffusion, while $\gamma>1$ corresponds to superdiffusive spreading \cite{metzler2000}. The same exponent can arise from different stochastic mechanisms; consequently, the mean square displacement alone does not determine the underlying process. 

\textit{3. Adsorption, desorption, absorption, and substrate.} 
Adsorption is the accumulation of molecules from a gas or liquid phase at a surface or interface. The accumulated species is the adsorbate, and the solid or interface on which accumulation occurs is the adsorbent or sorbent. Desorption is the reverse process. Adsorption should be distinguished from absorption: in absorption the species penetrates into the bulk of another material, whereas adsorption is primarily an interfacial phenomenon.

The adsorbed amount will be denoted by $q(t)$, with equilibrium value $q_e$. Classical adsorption kinetics is often described by pseudo-first-order or pseudo-second-order equations. In these phenomenological laws the approach to equilibrium is proportional, respectively, to the first or second power of the remaining adsorption capacity $q_e-q(t)$. The adjective \emph{pseudo} indicates that the fitted kinetic order need not coincide with the molecularity of an elementary chemical reaction.

In enzyme chemistry, the term \emph{substrate} denotes the molecule on which the enzyme acts. In electrochemistry we shall also use $S$ for an initial electroactive species when discussing a reaction such as $S\to P+ne^-$. The notation is contextual, but the physical role will always be stated explicitly.

\textit{4. Electrodes and electrochemical measurements.}
An electrode is an electronic conductor in contact with an ionic conductor (the electrolyte solution); electrochemical reactions occur at the electrode/electrolyte interface. A simple oxidation reaction can be written as $S\longrightarrow P+ne^-$, where $n$ electrons are transferred. The measured electric current is proportional to the molecular flux at the electrode. Up to a sign convention, Faraday's law gives
\[
I(t)=nFAJ(0,t),
\]
where $F$ is Faraday's constant and $A$ is the electrode area. This direct relation between a measurable electrical signal and a transport flux is one reason electrochemistry became an early and particularly natural field for fractional-calculus methods \cite{oldham2010}.

\textit{5. Surfactants and interfaces.}
A surfactant tends to accumulate at the boundary between two phases. Its adsorption changes the interfacial tension and can influence the dynamics of the entire multiphase system. In many experiments a molecule must first reach the interface by bulk transport and then cross an energetic or structural barrier before adsorbing. The observed kinetics therefore couples transport and interfacial reaction, and eliminating the bulk variable can generate a history dependence even if the original transport equation is classical.

\section{From local kinetics to memory and fractional operators}

Although numerous comprehensive monographs on fractional
calculus are available, we refer only to the classical
work of Podlubny \cite{podlubny1999} for the necessary
mathematical background, since the primary focus of this
review is on the application and mathematical justification
of fractional derivatives in chemical processes rather than
on the general theory of fractional calculus. General memory formulations, together with their fractional specializations, are surveyed from the viewpoint of non-local kinetics in \cite{hristov2023}. The simplest local relaxation law, $u'=-ku$, has an exponential solution and no explicit dependence on earlier states. A nonlocal description may instead involve a Volterra integral or a fractional derivative. The distinction between these representations is important because a singular power-law kernel can enter in more than one mathematically inequivalent way.

For $\beta>0$, define the Riemann--Liouville fractional integral by (see, e.g., \cite{podlubny1999})
\[
\Ii{\beta}f(t)=\frac{1}{\Gamma(\beta)}\int_0^t(t-s)^{\beta-1}f(s)\,ds.
\]
For $0<\alpha<1$ and $u\in AC([0,T])$, the Caputo derivative is (see, e.g., \cite{podlubny1999})
\begin{equation}\label{eq:caputo}
\CD{\alpha}u(t)=\Ii{1-\alpha}u'(t)=\frac{1}{\Gamma(1-\alpha)}\int_0^t(t-s)^{-\alpha}u'(s)\,ds .
\end{equation}
The Riemann--Liouville derivative of order $1-\alpha$ is $\RL{1-\alpha}u(t)=\frac{d}{dt}\Ii{\alpha}u(t)$ whenever the right-hand side is well defined (for example, in $L^1_{\rm loc}(0,T)$) (see, e.g., \cite{podlubny1999}). We use the usual causal convention for all fractional integrals and derivatives.

For $\lambda>0$, the fractional relaxation equation $\CD{\alpha}u(t)=-\lambda u(t)$, $u(0)=u_0$, is equivalent to
\begin{equation}\label{eq:volterra}
u(t)=u_0-\frac{\lambda}{\Gamma(\alpha)}\int_0^t(t-s)^{\alpha-1}u(s)\,ds .
\end{equation}
Thus the ordinary Volterra formulation contains the locally integrable memory kernel $t^{\alpha-1}$. Equivalently, $u'(t)=-\lambda\,\RL{1-\alpha}u(t)$. If this last identity is written only formally as $u'=-K*u$, then $K(t)$ is proportional to $t^{\alpha-2}$, which is not locally integrable for $0<\alpha<1$ and must be interpreted in the Riemann--Liouville or distributional sense. An ordinary convolution law with $K(t)\sim t^{-\alpha}$ should therefore not be identified directly with the Caputo relaxation equation.

The two-parameter Mittag--Leffler function is (see, e.g., \cite{gorenflo2020})
\begin{equation}\label{eq:ml}
E_{\alpha,\beta}(z)=\sum_{m=0}^\infty\frac{z^m}{\Gamma(\alpha m+\beta)},\qquad\alpha>0,\ \beta>0,
\end{equation}
and $E_\alpha(z)=E_{\alpha,1}(z)$ \cite{gorenflo2020}. Hence $u(t)=u_0E_\alpha(-\lambda t^\alpha)$. For $\alpha=1$ this is the exponential law. For $0<\alpha<1$,
\[
E_\alpha(-\lambda t^\alpha)\sim\frac{1}{\lambda\Gamma(1-\alpha)}t^{-\alpha},\qquad t\to\infty,
\]
so the long-time decay is algebraic. Pollard proved that $E_\alpha(-x)$ is completely monotone for $0<\alpha\le1$ and $x\ge0$ \cite{pollard1948}. 

These formulas suggest a basic modeling principle: the operator should represent a specified memory or transport mechanism, not merely an additional fitting parameter. The same mathematical operator may have different physical meanings depending on whether it arises from a waiting-time distribution, elimination of hidden variables, a constitutive flux law, or a distribution of ordinary relaxation modes.

A second distinction concerns the definition of flux. A time-fractional concentration equation can be paired either with a gradient-based Fick-type flux or with a memory flux derived from a local conservation law. These choices agree in the classical limit but can produce different boundary observables when $\alpha\neq1$. This issue becomes decisive in electrochemistry, where the measured current is itself a flux.

\section{Electrochemistry: from semi-integration to anomalous impedance}

\textit{1. Semi-infinite diffusion and the natural appearance of a half-order derivative.} Electrochemistry provides one of the clearest examples in which a fractional operator appears without fractionalizing the original bulk equation. Consider an electroactive species $S$ diffusing toward a planar electrode in a quiescent solution containing an excess of supporting electrolyte, so that migration and convection can be neglected. If the electrode radius of curvature is much larger than the diffusion-layer thickness, and the solution is much thicker than the diffusion layer, a one-dimensional semi-infinite model is appropriate:
\begin{equation}\label{eq:semi}
\frac{\partial C}{\partial t}=D\frac{\partial^2C}{\partial x^2},\qquad x>0.
\end{equation}
The coordinate $x=0$ represents the electrode surface. The concentration $C(x,t)$ is the amount of electroactive species per unit volume, while $D$ has units of length squared per unit time and controls the rate at which diffusion spreads concentration disturbances.

The surface flux is $J(0,t)=-D\,\partial_xC(0,t)$, and Faraday's law gives $I(t)=nFAJ(0,t)$. Thus the measured current is, up to known multiplicative factors and a sign convention, a direct measurement of the molecular flux at the boundary. This point is essential: the quantity that an electrochemical experiment records is not the concentration gradient by itself but the constitutively defined flux.

The precise boundary relation depends on the initial and far-field data. Assume
\[
C(x,0)=C_b\quad (x>0),\qquad \lim_{x\to\infty}C(x,t)=C_b,
\]
and prescribe the surface trace $C_s(t)=C(0,t)$. If $C_s-C_b$ is causal and of sufficient regularity for the Laplace transform (for instance, locally absolutely continuous and of exponential order), elimination of the bulk variable gives
\begin{equation}\label{eq:dtn}
\frac{\partial C}{\partial x}(0,t)
=-\frac{1}{\sqrt D}\,\RL{1/2}\bigl(C_s-C_b\bigr)(t)
=\frac{C_b}{\sqrt{\pi Dt}}-\frac{1}{\sqrt D}\,\RL{1/2}C_s(t).
\end{equation}
If, in addition, $C_s$ is continuous at $t=0$, $C_s(0)=C_b$, and $C_s\in AC([0,T])$, then
\[
\RL{1/2}C_s(t)=\CD{1/2}C_s(t)+\frac{C_b}{\sqrt{\pi t}},
\]
and \eqref{eq:dtn} reduces to
\begin{equation}\label{eq:dtn-caputo}
\frac{\partial C}{\partial x}(0,t)=-\frac{1}{\sqrt D}\,\CD{1/2}C_s(t).
\end{equation}
The continuity condition is not a technicality. In the ideal Cottrell experiment the electrode is switched instantaneously to a perfectly absorbing state, so $C_s(t)=0$ for $t>0$ while the bulk initial value is $C_b$; the boundary trace has a jump at $t=0$. Then the Caputo formula \eqref{eq:dtn-caputo} cannot be used literally with the post-step trace, whereas the Riemann--Liouville form \eqref{eq:dtn} remains the appropriate causal Dirichlet-to-Neumann relation and gives $\partial_xC(0,t)=C_b/\sqrt{\pi Dt}$. Thus the type of derivative and the compatibility of the initial and boundary data must be stated together. Consequently, the measured current and the surface concentration are related by a temporal half-order boundary operator, while the semi-integral of the current is linearly related to the concentration at the electrode.

Oldham and Spanier recognized this half-order structure already in the early development of semi-integral electroanalysis \cite{oldham1970}. Subsequent work developed the theory and experimental methodology of semi-integration for current and potential-sweep data \cite{grenness1972,goto1973,oldham1973,oldham1976,goto1977,imbeaux1973}. Oldham later emphasized the broader fractional-differential interpretation \cite{oldham2010}. More recently, semi-integral electroanalysis has been used to determine diffusion coefficients in highly viscous ionic liquids \cite{bentley2013}.

The mathematical meaning is especially transparent. The half-order operator is the boundary operator produced after eliminating the spatial variable from a classical, well-posed diffusion equation:
\[
\text{classical bulk PDE}\longrightarrow\text{elimination of }x\longrightarrow\text{half-order boundary operator}.
\]
The order $1/2$ is not fitted. It is fixed by the parabolic operator and the semi-infinite geometry. This mechanism provides a useful benchmark when judging more phenomenological fractional models.

\textit{2. Anomalous diffusion and fractional voltammetry.} When the transport itself is anomalous, fractional calculus may enter the bulk equation. Mirvceski and Tomovski considered cyclic voltammetry under a time-fractional diffusion law \cite{mirceski2009}. Their formulation uses a Riemann--Liouville-type time-fractional concentration equation in a semi-infinite planar geometry and imposes the current through the ordinary concentration gradient at the electrode,
\[
\frac{I(t)}{nFA D_f}=\left.\frac{\partial c_{\rm Red}}{\partial x}\right|_{x=0}
=-\left.\frac{\partial c_{\rm Ox}}{\partial x}\right|_{x=0},
\]
with homogeneous initial bulk concentrations, equal generalized diffusivities for the redox pair, and a reversible Nernst boundary condition during a linear potential sweep. Under this \emph{gradient-current closure}, the Laplace-transform solution expressed through Wright functions yields the Cottrell-type and voltammetric scalings
\begin{equation}\label{eq:mt}
I(t)\sim t^{-\alpha/2},\qquad I_p\sim v^{\alpha/2},
\end{equation}
for $0<\alpha\le1$. The first law refers to the corresponding semi-infinite step experiment, while the second is the peak-current law for the reversible linear sweep analyzed in \cite{mirceski2009}. For normal diffusion, $\alpha=1$, both reduce to the classical exponents $-1/2$ and $1/2$.

It is important that the gradient-current relation above is a constitutive choice of that model; it is not, by itself, the flux obtained by combining the same fractional concentration equation with the local conservation law. In a continuous-time-random-walk formulation \cite{metzler2000}, one instead preserves
\begin{equation}\label{eq:memflux}
\frac{\partial C}{\partial t}=-\nabla\cdot J_{\rm mem},\qquad
J_{\rm mem}=-d_\alpha\,\RL{1-\alpha}\nabla C .
\end{equation}
For standard initial data this formulation is equivalent to the Caputo diffusion equation $\CD{\alpha}C=d_\alpha\Delta C$, but the experimentally measured current must now be identified with the physical memory flux,
\[
I(t)=nFAJ_{\rm mem}(0,t),
\]
not with $-d_\alpha\partial_xC(0,t)$ alone. For a uniform initial concentration in a semi-infinite domain and a perfectly absorbing planar boundary after a potential step, $\partial_xC(0,t)\sim t^{-\alpha/2}$ and hence
\begin{equation}\label{eq:memscal}
I(t)\sim t^{\alpha/2-1},\qquad Q(t)=\int_0^t I(s)\,ds\sim t^{\alpha/2}.
\end{equation}
For a linear potential sweep with otherwise analogous diffusion-controlled assumptions, dimensional scaling gives $I_p\sim v^{1-\alpha/2}$. This sweep-rate exponent is therefore a prediction of the specified conservation law and current closure, not a universal consequence of the fractional diffusion equation alone.

The two closures coincide at $\alpha=1$ but lead to opposite departures from the classical slope $1/2$ for subdiffusion. In the gradient-current model, $d\log I_p/d\log v=\alpha/2<1/2$; in the conservative memory-flux model it is $1-\alpha/2>1/2$. Consequently, an experimental power-law slope cannot be converted into a fractional order until the geometry, initial data, electrode kinetics, potential protocol, and---most importantly---the constitutive definition of the measured flux have been fixed. 

\textit{3. Anomalous diffusion and electrochemical impedance.} Classical semi-infinite diffusion produces the Warburg behavior $Z(\omega)\propto(i\omega)^{-1/2}$. Bisquert and Compte derived anomalous-diffusion generalizations in which
\begin{equation}\label{eq:bc}
Z(\omega)\propto(i\omega)^{-\beta/2},
\end{equation}
and related $\beta$ to the transport law $\langle r^2(t)\rangle\sim t^\beta$ \cite{bisquert2001}. Their analysis is valuable because the impedance power law is connected to a model of anomalous transport rather than introduced as an isolated frequency-domain fitting element. The closure issue of Section~4.2 applies here as well. Under a semi-infinite linearization in which the concentration perturbation satisfies a time-fractional diffusion equation, the boundary gradient scales as $(i\omega)^{\beta/2}$. If the measured current is proportional to that gradient, one obtains \eqref{eq:bc}. If, instead, the current is the conservative memory flux in \eqref{eq:memflux}, the additional factor $(i\omega)^{1-\beta}$ changes the diffusional impedance scaling to
\[
Z_{\rm diff}(\omega)\propto(i\omega)^{-(1-\beta/2)}.
\]
Thus even in the same planar, linear, semi-infinite geometry the mapping from an impedance angle to a transport exponent depends on the current closure. Neither exponent should be transferred unchanged to finite-length diffusion, nonlinear electrode kinetics, or mixed interfacial control without a separate derivation.

The Poisson--Nernst--Planck framework has been extended further by introducing distributed-order time derivatives in the bulk and fractional operators in the boundary conditions \cite{lenzi2017}. In the low-frequency regime these models generate responses of the form $Z(\omega)\propto(i\omega)^{-\gamma}$. An important mathematical consequence is nonuniqueness of mechanism: a similar low-frequency slope may arise from memory in the bulk, memory at the interface, or a combination of both. Silva and coauthors obtained related anomalous responses for electrolytic cells containing several ion groups and distributed-order boundary effects \cite{silva2015}.

Bounded diffusion creates a different type of nonclassical frequency response. In finite electrode particles the response crosses from Warburg-like behavior at high frequency to a finite-length regime at low frequency. Gabano, Poinot, and Huard constructed fractional gray-box models for bounded diffusion in battery electrodes and used them to represent planar, cylindrical, and spherical geometries \cite{gabano2017}. In this setting the fractional element is an efficient reduced representation of a finite-geometry diffusion operator; its order should therefore not automatically be interpreted as a CTRW exponent.

\textit{4. Constant-phase elements and distributions of relaxation times.} A constant-phase element is commonly written as
\begin{equation}\label{eq:cpe}
Z_{\rm CPE}(\omega)=\frac{1}{Q_{\rm CPE}(i\omega)^\nu},\qquad0<\nu\le1 .
\end{equation}
where $Q_{\rm CPE}$ is the CPE amplitude. This notation is kept distinct from the transferred charge $Q(t)$ used in \eqref{eq:memscal}. For $\nu=1$ the response is ideally capacitive; for $\nu=1/2$ it has the form of the Warburg impedance of semi-infinite diffusion, which is precisely the frequency-domain counterpart of the half-order boundary operator \eqref{eq:dtn}; for other $0<\nu<1$ the phase is nearly constant over a broad frequency interval. Balaguera analyzes the relation among CPE parameters, fractional calculus, and distributions of relaxation times \cite{balaguera2025}. The same nonideal capacitive response can often be represented either by a fractional constitutive law or by a positive distribution of ordinary relaxation modes.

\begin{figure}[t]
\centering
\includegraphics[width=\textwidth]{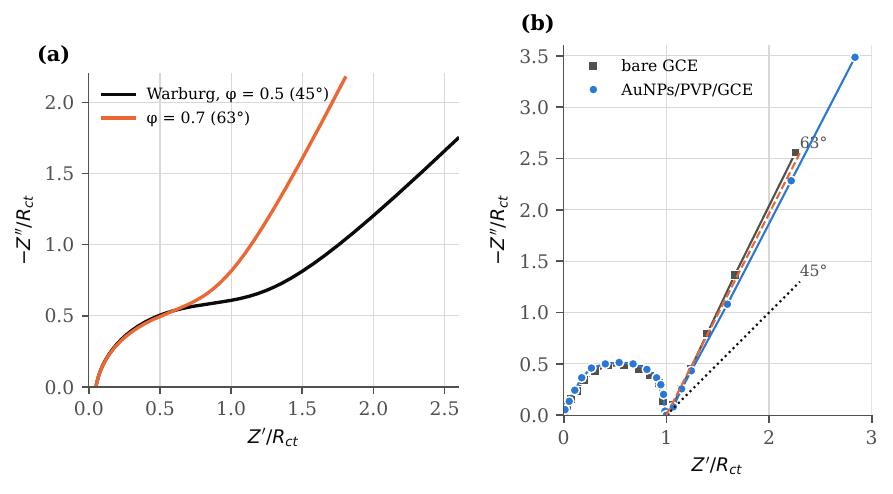}
\caption{(a) Simulated Nyquist plots of a Randles circuit with a CPE double layer ($\nu=0.9$) and a low-frequency element $Z_d=\sigma(i\omega)^{-\varphi}$: the ideal Warburg case $\varphi=0.5$ and $\varphi=0.7$. (b) Spectra of bare GCE and AuNPs/PVP/GCE in 5~mM [Fe(CN)$_6$]$^{3-/4-}$/0.1~M KCl at the open-circuit potential, re-plotted from \cite{rajabova2026} (digitised), each normalised by its own charge-transfer resistance. Both axes have equal scales; guide lines mark $45^\circ$ and $63^\circ$.}
\label{fig:eis}
\end{figure}

This equivalence shifts the mathematical focus from existence of an initial-value problem to causality, passivity, complete monotonicity, and identifiability of the relaxation spectrum. A fractional fit is physically more persuasive when the corresponding time-domain response is compatible with a positive relaxation measure. Conversely, if several distinct positive spectra or circuit topologies fit the same frequency window, the fitted order should be interpreted as an effective descriptor rather than a uniquely identified microscopic quantity. Figure~\ref{fig:eis}(a) shows how a noninteger low-frequency exponent changes a simulated Randles spectrum; panel (b) shows experimental spectra discussed in Section~\ref{sec:case}.

\section{Adsorption, desorption, and sorption with memory}

\textit{1. From pseudo-first-order kinetics to Mittag--Leffler relaxation.} A standard pseudo-first-order adsorption law is
\begin{equation}\label{eq:pfo}
\frac{dq}{dt}=k(q_e-q),
\end{equation}
with solution $q(t)=q_e-(q_e-q_0)e^{-kt}$. If the measured approach to equilibrium has a much broader tail than a single exponential, one possible generalization is
\begin{equation}\label{eq:fpfo}
\CD{\alpha}q=k_\alpha(q_e-q),\qquad0<\alpha<1,
\end{equation}
which yields $q(t)=q_e-(q_e-q_0)E_\alpha(-k_\alpha t^\alpha)$.

Bakalis and Zerbetto compare classical, fractal, and fractional descriptions of adsorption kinetics \cite{bakalis2025}. Their distinction is conceptually useful: geometric or temporal heterogeneity may generate nonclassical kinetics without true memory, whereas a fractional derivative explicitly introduces a nonlocal dependence on the history of the process. For the linear equation \eqref{eq:fpfo}, existence and uniqueness follow directly from the explicit Mittag--Leffler representation; for nonlinear higher pseudo-orders, a recursion generating numerical values is not, by itself, a uniqueness theorem for the original nonlinear fractional equation.

The classical Ward--Tordai description of diffusion-controlled surfactant adsorption is already nonlocal because the amount adsorbed at the interface depends on the entire history of bulk diffusion. Hristov recast this structure in fractional form and showed how half-order operators arise in the Ward--Tordai setting \cite{hristov2016}. The analogy with semi-infinite electrochemistry is strong: a boundary history law emerges after solving a bulk diffusion problem rather than from a formal change of derivative in a surface ODE.

\textit{2. Heavy-metal sorption and statistical identification of the fractional order.} Friesen and coauthors modeled lead sorption by water hyacinth using a fractional relaxation model and, importantly, complemented the fit with parameter-uncertainty analysis \cite{friesen2015}. They calculated covariance information and joint confidence regions for the fractional order and kinetic parameter and also considered a setting in which one experimental data set was used for parameter estimation and another for validation.

This type of analysis addresses a weakness common to fractional fitting. If $\alpha$ and $k_\alpha$ are both free, a reduction in residual error does not show that $\alpha$ is statistically identifiable or genuinely separated from the integer value $1$. Confidence regions, profile likelihoods, correlation analysis, and out-of-sample prediction are needed to determine whether the fractional order behaves as a reproducible parameter rather than as additional curve-fitting flexibility.

\textit{3. SPR adsorption and nonexponential sensorgrams.} Surface plasmon resonance (SPR) measures binding at a sensor surface in real time. In the simplest one-to-one model, association and dissociation are governed by first-order ODEs and therefore produce exponential phases. Ferreira and coauthors applied a Caputo model to experimental SPR sensorgrams for the interaction of immobilized Baru protein with Congo Red dye \cite{ferreira2026}. In the comparisons reported there, the fractional one-to-one model reduced the fitting error by more than one order of magnitude relative to the corresponding integer-order one-to-one model and reproduced a rapid initial response followed by a much slower approach or decay. This is a quantitative improvement within that pair of models; it is not, by itself, evidence that alternative classical mechanisms such as mass-transfer limitation or surface heterogeneity have been excluded.

The experiment is particularly informative because the measured signal depends strongly on concentration and is clearly nonexponential. At the same time, SPR is sensitive not only to intrinsic binding but also to limitations of mass transfer, surface heterogeneity, and conformational rearrangements. A fractional one-to-one law can therefore act as an effective compression of several unresolved mechanisms. A decisive model comparison should include classical models with mass-transfer limitations, heterogeneous binding sites, and multiple states, not only the integer-order one-to-one model. Better fitting establishes usefulness; discrimination among mechanisms requires additional observables or experimental perturbations.

\textit{4. Fractional Langmuir dynamics and relaxation spectra.} Barbero, Evangelista, and Lenzi studied frequency dispersion in a fractional Langmuir approximation and compared fractional kinetics with descriptions based on distributions of relaxation times \cite{barbero2020}. In both representations effective adsorption and desorption parameters can become frequency dependent. The mathematical point is the same one encountered for CPEs: over a finite frequency interval, a fractional kernel may be nearly indistinguishable from a positive mixture of classical relaxation modes.

This creates an inverse problem that is more subtle than estimation of one exponent. One must ask whether the fractional representation is unique within a physically admissible class and whether the inferred kernel is compatible with complete monotonicity, positivity, and known limiting behavior. If a broad distribution of ordinary time constants gives an equivalent description, the fractional order should be viewed as a compact effective parameter unless independent data distinguish the two mechanisms.

\textit{5. Barrier--diffusion controlled adsorption under anomalous transport.} Bazhlekov and Bazhlekova developed a model in which anomalous behavior is introduced at the level of the constitutive flux and then propagated consistently to the interface \cite{bazhlekov2025,bazhlekov2026}. In each bulk phase,
\begin{equation}\label{eq:bb1}
\partial_tC^\pm=-\nabla\cdot J^\pm,\qquad
J^\pm=-d^\pm\,\RL{1-\gamma^\pm}\nabla C^\pm,
\end{equation}
with $0<\gamma^\pm<2$. For $0<\gamma<1$, $\RL{1-\gamma}$ is a Riemann--Liouville derivative; for $1<\gamma<2$, the same notation means the fractional integral $\Ii{\gamma-1}$. The classical Fick law is recovered at $\gamma=1$.

The one-dimensional semi-infinite formulation in \cite{bazhlekov2026} is supplied with uniform bulk initial data and far-field conditions,
\[
C^\pm(x,0)=C_e^\pm,\qquad
\lim_{x\to\pm\infty}C^\pm(x,t)=C_e^\pm,
\]
together with interface traces $C_s^\pm(t)=C^\pm(0^\pm,t)$. The accumulated interfacial amounts satisfy
\begin{equation}\label{eq:bb2}
J_s^\pm=K_a^\pm\Bigl(C_s^\pm(t)-\frac{f(\sigma(t))}{K^\pm}\Bigr),\qquad
\frac{d\sigma^\pm}{dt}=J_s^\pm,
\qquad
\sigma(t)=\sigma^+(t)+\sigma^-(t),
\end{equation}
with $\sigma^\pm(0)=0$. Here $K_a^\pm>0$ are adsorption-rate parameters, $K^\pm>0$ are equilibrium parameters, and $f$ specifies the adsorption isotherm.

The range $1<\gamma<2$ deserves an explicit comment. Combining \eqref{eq:bb1} with the local balance and taking Laplace transforms gives
\[
s\widehat C-C_e=d\,s^{1-\gamma}\Delta\widehat C.
\]
After multiplication by $s^{\gamma-1}$, this is the Laplace-domain form of the Caputo equation $\CD{\gamma}C=d\Delta C$ with
\[
C(x,0)=C_e,\qquad \partial_tC(x,0)=0,
\]
when $1<\gamma<2$. Thus the zero initial rate is already encoded when the first-order balance law with the integral constitutive flux is rewritten as a differential equation of order $\gamma$; it is not an independent extra physical law. For $0<\gamma<1$, the same calculation gives the usual Caputo equation of order $\gamma$ with only the initial concentration prescribed.

Elimination of the semi-infinite bulk variables yields
\begin{equation}\label{eq:bb3}
C_s^\pm(t)=C_e^\pm-\frac{1}{\sqrt{d^\pm}}\,\RL{\gamma^\pm/2}\sigma^\pm(t),
\end{equation}
where $0<\gamma^\pm/2<1$. Substitution into \eqref{eq:bb2} gives
\begin{equation}\label{eq:bb4}
\frac{d\sigma^\pm}{dt}
+\frac{K_a^\pm}{\sqrt{d^\pm}}\,\RL{\gamma^\pm/2}\sigma^\pm
=K_a^\pm C_e^\pm-\frac{K_a^\pm}{K^\pm}f(\sigma),
\qquad \sigma^\pm(0)=0.
\end{equation}
The derivation in \cite{bazhlekov2026} is carried out in a Laplace-transform setting. A natural time-domain class for \eqref{eq:bb4} is $\sigma^\pm\in AC_{\rm loc}([0,T])$ with $\RL{\gamma^\pm/2}\sigma^\pm\in L^1_{\rm loc}(0,T)$; in this class the zero initial values also remove the usual Riemann--Liouville singular initial terms. A complete nonlinear well-posedness theorem should state such a solution class explicitly.

The fractional boundary operator is therefore derived from anomalous bulk transport rather than inserted directly into the adsorption rate law:
\begin{align*}
&\text{anomalous transport}\to\text{fractional flux}\to\text{fractional PDE}\\
&\quad\to\text{elimination of the bulk variable}\to\text{fractional adsorption equation}.
\end{align*}
For a Henry isotherm the system is linear and admits explicit representations in generalized Mittag--Leffler functions. Volmer and van der Waals isotherms lead to nonlinear Volterra equations treated numerically by predictor--corrector procedures. The cited analysis does not state a general global existence-and-uniqueness theorem for all nonlinear isotherms in a function class sufficient to justify every numerical trajectory. This gap motivates Open Problems~1 and~2 below.

Lenzi and coauthors considered another sorption problem governed by fractional diffusion and obtained Green-function and mean-square-displacement results \cite{lenzi2016}. Together, the adsorption literature displays at least three distinct origins of memory: anomalous bulk transport, a spectrum of interfacial relaxation times, and boundary history generated by elimination of spatial diffusion. Similar kinetic curves can therefore correspond to different mathematical mechanisms.

\section{Reaction--subdiffusion and chemical transformations in porous media}

\textit{1. Methanol transport in mesoporous materials: geometry versus adsorption.} A particularly instructive sequence of experimental and modeling studies concerns methanol transport in mesoporous silica and zeolite-containing catalyst pellets \cite{zhokh2017a,zhokh2017b,zhokh2017c,zhokh2017d,zhokh2019}. In the cited studies, classical Fickian diffusion was reported to be insufficient for several transient transport curves, while time-fractional descriptions gave lower discrepancies for the measured non-Fickian kinetics within the model comparisons performed there. The studies are valuable because the fractional exponent is not examined in isolation: its dependence on pore size, material composition, and adsorption properties is compared across related systems.

The interpretation is more subtle than a simple identification of the exponent with geometric fractality. Changes in adsorption energy and surface chemistry can influence the inferred anomalous exponent at least as strongly as substantial changes in pore size. Trapping and desorption therefore compete with geometry in determining the effective transport law. This observation is chemically important: in a porous catalyst the waiting-time statistics of a molecule are determined not only by how it moves through the pore network but also by how long it remains bound to the internal surface.

A later comparison with fractional Brownian motion makes the issue sharper \cite{zhokh2019}. A time-fractional CTRW model and fractional Brownian motion can share the same mean-square-displacement law while possessing different stochastic structures. The former is naturally associated with heavy-tailed waiting times and aging; the latter has correlated, Gaussian increments. Thus the exponent in $\langle x^2(t)\rangle\sim t^\alpha$ does not identify the stochastic model. First-passage statistics, aging tests, response to changes in adsorption energy, or trajectory-level information are needed to distinguish mechanisms.

\textit{2. Thiele modulus and reactor-level consequences.} For a porous catalyst particle, the classical Thiele modulus compares the characteristic rate of reaction with the rate of diffusion inside the particle. Its importance is not merely mathematical: it controls the effectiveness factor, namely the fraction of the catalyst volume that contributes effectively to reaction. Zhokh and Strizhak generalized this engineering concept to anomalous transport \cite{zhokh2018}. Once the transport time no longer scales quadratically with length in the classical way, the balance between reaction and intraparticle transport changes and the usual dimensionless criterion must be modified.

This type of result goes beyond fitting a transient concentration curve. A fractional transport model can alter predicted catalyst effectiveness, optimal particle size, and the regime in which internal transport becomes rate limiting. Such reactor-level quantities provide a stronger route to experimental validation because they can be tested independently of the transient data used to estimate the anomalous exponent.

Zhang, Li, Xiu, and Rodrigues considered Riemann--Liouville diffusion--reaction equations in fractal porous media and derived analytical and numerical solutions for first-order irreversible reactions with intraparticle diffusion and external mass-transfer resistance \cite{zhang2021}. Their formulas for effectiveness factors illustrate how fractional transport enters directly into a classical chemical-engineering quantity. The combination of anomalous diffusion, external resistance, and reaction also shows why dimensional analysis is essential: a generalized diffusion coefficient, a reaction-rate parameter, and a fractional time scale must be scaled consistently.

\textit{3. Why reaction and subdiffusion cannot always be coupled by addition.} The most important conceptual difficulty in reaction--subdiffusion is that chemical reaction can occur during the waiting periods that generate subdiffusion. Transport and reaction are then statistically coupled. Seki, Wojcik, and Tachiya derived fractional reaction--diffusion equations from continuous-time random walks and showed that the reaction term may inherit the same memory structure as the transport term \cite{seki2003a,seki2003b}. In recombination problems this coupling is necessary to preserve the dependence of the reaction probability on intrinsic reactivity.

Consequently, an equation such as
\[
\CD{\alpha}C=D_\alpha\Delta C-kC
\]
should not be regarded as a universal reaction--subdiffusion law merely because it combines a fractional diffusion operator with a familiar reaction term. The correct form depends on what happens to the internal reaction clock while the particle is trapped. If the particle can react during a waiting period, if reaction terminates the waiting event, or if the reaction itself has a separate memory, different macroscopic equations result.

This point has been emphasized again in the recent survey of non-Fickian reaction--diffusion equations by Zhokh and Strizhak \cite{zhokh2026}. Their comparison of CTRW, L\'evy-flight, Cattaneo, fractal, and Feynman--Kac formulations shows that the apparently simple operation of ``adding a reaction term'' is physically safe only in specific classes of non-Fickian transport. In memory-dominated subdiffusion, an ad hoc coupling may violate positivity or produce reaction kinetics inconsistent with the underlying stochastic process.

Traytak's fractional differentiation method gives another perspective on the same class of problems \cite{traytak2024}. For several subdiffusion-controlled trapping reactions, reaction-rate coefficients can be obtained by a factorization procedure without solving the complete initial-boundary-value problem first. Here fractional differentiation is not only a modeling language but also an analytic tool for extracting an observable quantity directly from a transport problem.

\textit{4. Green functions, short and long times.} Fractional reaction--diffusion equations often possess qualitatively different short- and long-time regimes because the Mittag--Leffler function crosses from near-exponential behavior to an algebraic tail. Zhokh and Strizhak derived Green-function approximations on different time scales and identified a crossover time for a linear time-fractional reaction--diffusion equation \cite{zhokh2023}. They applied the resulting formulas to a heterogeneous catalytic conversion problem.

This scale dependence has an important inverse-problem consequence. A fractional order fitted over a narrow early-time window need not coincide with the effective order inferred from a long-time tail if the model itself contains crossovers or if the data mix several transport mechanisms. Reported fractional parameters should therefore be accompanied by the observation window and by sensitivity tests with respect to that window.

\textit{5. Stoichiometric balance in nonlinear fractional kinetics.} Fractional reaction models must preserve the stoichiometric structure of the underlying chemistry. In a classical reaction network, $\dot c=S\,r(c)$, where $S$ is the stoichiometric matrix and $r(c)$ is the vector of reaction rates. Any vector $\ell$ satisfying $\ell^TS=0$ generates a linear invariant $\ell^Tc(t)=\mathrm{const}$ in the absence of external sources and sinks.

Barajas-Villarruel, Dami\'an-V\'azquez, and Rico-Ramirez proposed a general route to higher-order fractional kinetics designed to retain mass balance \cite{barajas2025}. Rather than fractionalizing each species equation independently, the construction modifies the rate description while keeping the same stoichiometric reaction flux in all component balances. Reversible second-order kinetics and heterogeneous catalytic reactions are used as examples.

This structural viewpoint is preferable to component-wise fractionalization. If different species are assigned unrelated fractional operators, as in several fractional versions of the Brusselator systems \cite{asiedu2025}, a classical conservation law may disappear even though the reaction mechanism itself has not created or destroyed the conserved atomic component. A satisfactory fractional reaction-network theory should therefore specify simultaneously the memory operator, the stoichiometric architecture, positivity of concentrations, and the class of admissible rate laws.

\section{Fractional chemical kinetics and reaction networks}

\textit{1. Linear and sequential reactions.} The simplest fractional kinetic law replaces exponential relaxation by a Mittag--Leffler law. For a sequence of first-order reactions, however, memory propagates through intermediate species and changes not only decay rates but also the location and height of concentration maxima. Murga, Pskhu, and Rekhviashvili studied sequential reactions in a fractional integro-differential setting and obtained exact expressions involving Mittag--Leffler and Prabhakar functions \cite{murga2026}. Such formulas provide useful benchmarks for numerical algorithms and clarify how a nonlocal kernel modifies the transfer of mass through a reaction chain.

A separate question is the physical origin of the memory. Qureshi and Aziz investigated a three-species batch-reactor system with Caputo derivatives and obtained closed-form concentration profiles by the Laplace transform \cite{qureshi2020}. The analysis carefully tracks the dimensions of the fractional parameters and shows how the trajectories vary with the order. Mathematically, this is a transparent solvable example. Chemically, the dependence of curves on $\alpha$ demonstrates sensitivity to memory but does not by itself identify a microscopic mechanism; experimental calibration would still be required.

\textit{2. Photocatalytic degradation.} Photocatalytic degradation is naturally multiscale. Light absorption, generation of reactive intermediates, transport to the catalyst surface, adsorption, reaction, and desorption may all contribute to the measured disappearance curve. Wang compared first-order, second-order, and fractional first-order descriptions for experimental degradation data and reported lower fitting discrepancies for the fractional relaxation law on several data sets \cite{wang2018}. The order and a characteristic degradation time were then used as material descriptors.

A better fit is suggestive but not unique evidence for fractional memory. Nonexponential degradation can also be generated by changing light intensity inside the sample, heterogeneous active sites, catalyst deactivation, multi-step chemistry, or mass-transfer limitations. A mechanistically discriminating experiment would vary one of these factors while testing whether a single fractional order remains transferable. If the same $\alpha$ predicts different initial concentrations, catalyst loadings, irradiation protocols, and hydrodynamic conditions, its physical interpretation becomes much stronger.

\textit{3. Numerical methods for fractional reaction equations.} A large body of recent work develops numerical or semi-analytical solvers for fractional equations motivated by chemical systems. Sylvia and Ghosh proposed a Haar-wavelet method for a nonlinear three-variable fractional chemical-reaction model and analyzed convergence, stability, and existence/uniqueness of the stated problem \cite{sylvia2024}; they later applied a related wavelet framework to a chemical-clock model \cite{sylvia2025}. Manohara and Kumbinarasaiah used Fibonacci-wavelet collocation for a fractional Brusselator system \cite{manohara2024}.

For reaction--diffusion equations, Chawla, Kumar, and Vigo-Aguiar developed a method for a fourth-order multi-term time-fractional problem with weak initial singularities and analyzed stability and convergence on graded meshes \cite{chawla2025}. Roul and Rohil constructed a fourth-order compact ADI scheme for a two-dimensional time-fractional reaction--subdiffusion equation, with fourth-order spatial accuracy and a time order governed by the fractional singularity \cite{roul2024}. Nayak and Jena considered fractional Fokker--Planck equations motivated by stochastic molecular transport and compared computational errors and convergence indicators for transform/decomposition techniques \cite{nayak2026}. Fractional optimal-control problems for coupled diffusion systems have likewise been treated by transform-based methods \cite{sebaq2025}.

These contributions are important for computation, but the type of conclusion must be stated precisely. A convergence theorem shows that a numerical approximation approaches the solution of a specified mathematical model. It does not establish that the model is chemically correct. Conversely, experimental agreement cannot replace convergence analysis. The strongest computational studies connect both levels: they begin with a physically justified model, establish enough well-posedness to identify the target solution, and then prove that the discretization converges to that solution.

\section{Oscillatory reactions and complex dynamics}

Oscillatory reactions provide a natural setting in which memory can affect stability thresholds, transient amplitudes, phase relations, and the response to forcing. They also illustrate a danger of fractional modeling: adding a free order $\alpha$ to an already nonlinear system can generate rich dynamics even when no independent memory mechanism has been identified.

\textit{1. Continuous-time fractional models.}
Several works treat Brusselator- or Belousov--Zhabotinsky-type systems as continuous-time equations with Caputo derivatives. Asiedu and coauthors applied decomposition and variational-iteration techniques to fractional Brusselator equations \cite{asiedu2025}, while Manohara and Kumbinarasaiah used wavelet collocation \cite{manohara2024}. For the BZ system, Alshehry and coauthors combined a Caputo model with an optimal auxiliary-function method \cite{alshehry2023}; Rehman and coauthors compared transform/asymptotic approaches for a fractional Noyes--Field system \cite{rehman2024}; and AlBaidani and Alzahrani proposed further approximations for a time-fractional BZ model \cite{albaidani2026}. These studies show how prescribed fractional orders alter computed oscillatory profiles and stability properties. Such conclusions concern the chosen continuous fractional equations; they are not, by themselves, experimental evidence for power-law memory in the underlying chemistry.

\textit{2. A discrete-time fractional Brusselator map.}
Berkal and Almatrafi study a different mathematical object \cite{berkal2025}. They begin with a Caputo-time modified Brusselator, discretize it by a piecewise-constant-argument construction, and then analyze the resulting discrete map. The reported period-doubling and Neimark--Sacker bifurcations, invariant closed curves, Lyapunov exponents, and $0$--$1$ chaos test therefore belong to the \emph{discrete-time system}. In particular, a Neimark--Sacker bifurcation is the discrete analogue of a Hopf bifurcation and should not be described as a bifurcation theorem for the original continuous Caputo differential system. The distinction matters because stability regions and bifurcation mechanisms for fractional differential equations and for their discrete maps are not interchangeable.

The unresolved chemical question is common to both settings: is a fitted or prescribed fractional order an experimentally reproducible characteristic of the chemical mechanism, or does it merely shift nonlinear regimes of a chosen model? A convincing identification would require the same order to predict several independent observables---for example oscillation period, amplitude relaxation, phase lag under forcing, and, where applicable, spatial wave speed---over changes in temperature and initial composition. Chemical-clock models raise the same issue. Scaling across families of experiments is more discriminating than reproduction of a single numerical trajectory.

\section{Enzyme and biochemical kinetics}

\textit{1. Why memory may arise in enzyme kinetics.} 
Classical Michaelis--Menten kinetics is based on
\[
E+S\rightleftharpoons ES\longrightarrow E+P
\]
and, under standard quasi-steady assumptions, leads to a local rate law. Real enzymatic systems can contain conformational switching, molecular crowding, rebinding, heterogeneous active states, and broad residence-time distributions. These mechanisms can in principle generate memory, but the form of the memory kernel should be tied to a concrete mechanism rather than inferred solely from the improved flexibility of a fractional model.

Radhakrishnan, Chandru, and Tilahun considered fractional biochemical reaction models and used homotopy-based methods together with numerical and statistical comparisons of solutions for different prescribed orders \cite{radhakrishnan2025}. A subsequent study focused on a fractional Michaelis--Menten formulation \cite{radhakrishnan2026}. These works show how fractional orders change the substrate, enzyme--substrate-complex, and product profiles. Statistical differences among trajectories generated with different chosen values of $\alpha$, however, are not the same as statistical evidence from experimental data that $\alpha\neq1$. The latter requires estimation of the order, an uncertainty interval, and comparison with competing models.

Zulqarnain and coauthors developed a cooperative enzyme model with a Liouville--Caputo derivative, proved existence and uniqueness, and compared numerical procedures with a neural-network approximation trained on model-generated data \cite{zulqarnain2026}. The well-posedness analysis is a useful mathematical addition. The neural-network comparison, on the other hand, validates the approximation of the chosen equations rather than the biochemical mechanism itself. Experimental concentration or single-molecule data would be needed for mechanistic validation.

\textit{2. Variable-order models and delays.} Agilan and coauthors studied a variable-order Caputo enzyme model with time delays and analyzed existence, uniqueness, boundedness, stability, and numerical behavior \cite{agilan2025}. A time-dependent order $\delta(t)$ is intended to represent an evolving memory strength, while delays account for finite lags associated with intermediate or conformational processes. Because several inequivalent definitions of a variable-order derivative exist, the operator must be stated explicitly. One commonly used Type-I Caputo definition is
\begin{equation}\label{eq:vo-caputo}
{}^{\mathrm C}\!D_t^{\delta(t)}u(t)
=\frac{1}{\Gamma(1-\delta(t))}\int_0^t(t-s)^{-\delta(t)}u'(s)\,ds,
\qquad 0<\delta_-\le \delta(t)\le\delta_+<1,
\end{equation}
for $u\in AC([0,T])$. Other variable-order definitions place the time-dependent order inside the kernel in a different way and are not automatically equivalent to \eqref{eq:vo-caputo}.

Variable-order formulations create a serious inverse problem. If $\delta(t)$, several delay times, and kinetic constants are all unknown, their effects on the same concentration curve can be highly correlated. Moreover, a slowly varying order may mimic a slowly varying classical rate constant. Structural identifiability and independent measurements of characteristic delays are therefore essential before a reconstructed function $\delta(t)$ can be assigned a biochemical interpretation. Dimensional consistency also becomes nontrivial: since ${}^{\mathrm C}D_t^{\delta(t)}$ has the instantaneous dimension $T^{-\delta(t)}$, a coefficient multiplying the operator must be defined relative to a reference time scale or otherwise carry the corresponding variable dimension.

The broader numerical literature on variable-order equations also provides useful computational tools; examples include operational-matrix and hybrid-function methods applied to chemical models \cite{shiri2025}. Such flexibility is most meaningful when the order law is derived from an evolving physical environment or estimated under identifiable constraints. Otherwise an arbitrary function $\delta(t)$ can absorb substantial model discrepancy.

\section{Thermochemical, materials, and environmental applications}

\textit{1. Pyrolysis and thermal memory.} Pathak, Bhatt, and Chavada modeled non-isothermal biomass pyrolysis with a Caputo derivative \cite{pathak2026}. The model describes a multistage conversion of biomass into primary products and subsequent secondary transformations. The authors establish nonnegativity, boundedness, existence and uniqueness, and Ulam--Hyers-type stability before applying a fractional Runge--Kutta procedure and studying heating-rate effects.

This combination of analysis and process optimization is stronger than a purely graphical fractionalization. Nevertheless, the physical meaning of a ``thermal-memory'' order should be tested across several heating programs. If one value of $\alpha$ calibrated at one heating rate predicts conversion curves at other rates, the model has genuine predictive content. If the order must be refitted independently for every temperature program, it functions mainly as an effective surrogate parameter.

\textit{2. Alkali--silica reactions and material aging.} Fractional models have also been proposed for alkali--silica reaction dynamics in cementitious materials \cite{kumar2022}. Such systems are intrinsically multiscale: chemical reaction, moisture transport, pore structure, swelling products, and evolving damage interact over long periods. A memory term is therefore plausible, but its physical interpretation is ambiguous. The observed hereditary effect may belong to chemical transformation, moisture redistribution, damage evolution, or an effective coarse-graining of all three.

This setting suggests a more systematic route to fractional models. Instead of fractionalizing an empirical reaction law directly, one can begin with a heterogeneous transport--reaction--damage system and ask whether homogenization or elimination of hidden microstructural variables produces a nonlocal-in-time effective equation. Such derivations would make the fitted order traceable to the material architecture.

\textit{3. Environmental chemistry and transport.} Naik and coauthors constructed a Caputo system coupling marine chemical variables, anthropogenic forcing, CO$_2$, and ocean pH and investigated boundedness, stability, and numerical behavior \cite{naik2025}. At the scale of regional environmental time series, however, long memory can result from seasonal forcing, spatial averaging, unresolved circulation, and slowly varying external drivers. A fitted fractional order therefore need not represent microscopic chemical memory. Robustness under changes in temporal resolution and forcing assumptions is a necessary diagnostic.

Bohaienko and coauthors introduced fractional-order derivatives into a vertically resolved model of carbon sequestration in coastal wetlands, combining memory with moisture dynamics and sediment layering \cite{bohaienko2025}. Such models can represent slow carbon turnover compactly, but a superposition of classical compartments with widely separated turnover times may produce similar observations. Model discrimination should therefore compare fractional kernels with positive multi-compartment alternatives.

\textit{4. A note on fractional quantum chemistry.} Fractional ideas extend beyond macroscopic transport and kinetics. Mostafanejad developed a perspective on fractional quantum chemistry that includes fractional Schr\"odinger equations, fractional Hartree--Fock and Kohn--Sham constructions, generalized molecular integrals, and variable-order operators intended to represent nonlocal electronic effects \cite{mostafanejad2021}. The role of the fractional operator is fundamentally different here: it modifies the quantum-mechanical operator structure rather than serving as a temporal memory kernel for a concentration.

This distinction matters conceptually. The mathematical questions concern operator domains, variational stability, correlation energies, numerical convergence, and physical calibration of nonlocal electronic operators. These problems belong to the broad landscape of fractional methods in chemistry, but their interpretation should not be conflated with anomalous diffusion or hereditary chemical kinetics.

\section{Dielectric relaxation, polymers, and distributed time scales}

Polymer relaxation provides a clear example of dynamics distributed over many molecular time scales. The Debye model assumes a single characteristic relaxation time. Real polymers contain segmental motions, interfacial polarization, polydispersity, and cooperative rearrangements, so their dielectric spectra are broad and often asymmetric. 

Renter\'ia-Balti\'errez and coauthors proposed an electrical fractional model based on cap-resistor elements and the complex electric modulus to describe both $\alpha$-relaxation and Maxwell--Wagner--Sillars interfacial polarization \cite{renteria2025}. For order $a$, $0<a<1$, the cap-resistor obeys $V(t)=(\tau_a^a/C)\,\CD{a}Q(t)$. The fitted fractional orders vary between relaxation mechanisms and are interpreted as measures of dispersity and molecular mobility. Because the model is fitted in both frequency and temperature domains, it provides more information than a single-frequency power law. 

The underlying mathematical issue is the representation of a broad spectrum by a small number of fractional parameters. Completely monotone relaxation functions can be represented as positive mixtures of exponentials, so a fractional law may act as a compact surrogate for a continuum of Debye processes. This is both a strength and a source of nonuniqueness: a Havriliak--Negami-type law, a fractional constitutive equation, and a finite or continuous mixture of exponentials can be almost indistinguishable on a limited frequency interval.

Accordingly, physically admissible fractional dielectric models should be checked for causality and passivity, and parameter identification should use broadband data. The same mathematical framework links polymer relaxation, CPE impedance, and fractional Langmuir kinetics through Bernstein functions, completely monotone kernels, and inverse Stieltjes or Laplace-transform problems.

\section{Non-Markovian atmospheric methane chemistry}

Chishtie proposed a fractional description of atmospheric methane motivated by multiscale transport, heterogeneous OH fields, heavy-tailed residence times, and continuous-time random walks \cite{chishtie2026}. The classical box balance
\[
\frac{dm}{dt}=E-\lambda m
\]
contains the methane burden $m(t)$, an effective emission rate $E$, and a loss rate $\lambda$. A power-law residence-time distribution motivates a replacement of the exponential relaxation by a Mittag--Leffler law and is supported in \cite{chishtie2026} by Monte Carlo subordination calculations.

There is, however, an important operator issue in the cited paper. The derivative in \cite{chishtie2026} is called Riemann--Liouville, but the subsequent Laplace-transform step is written as
\[
\mathcal L\{D_t^\alpha m\}(s)
=s^\alpha\widehat m(s)-s^{\alpha-1}m_0,
\]
which is the transform formula for the Caputo derivative (or, equivalently, for a regularized Riemann--Liouville derivative with the classical initial value built in), not for the standard Riemann--Liouville derivative. The Mittag--Leffler solution displayed there therefore corresponds mathematically to
\begin{equation}\label{eq:methane-caputo}
\CD{\alpha}m(t)=E-\lambda_\alpha m(t),\qquad m(0)=m_0,
\end{equation}
namely
\begin{equation}\label{eq:methane-sol}
m(t)=m_0E_\alpha(-\lambda_\alpha t^\alpha)
+\frac{E}{\lambda_\alpha}\bigl[1-E_\alpha(-\lambda_\alpha t^\alpha)\bigr],
\end{equation}
with $[\lambda_\alpha]=T^{-\alpha}$.

If one instead interprets the model with the standard Riemann--Liouville derivative,
\begin{equation}\label{eq:methane-rl}
\RL{\alpha}m(t)=E-\lambda_\alpha m(t),
\end{equation}
then the natural initial datum is fractional,
\[
\bigl(\Ii{1-\alpha}m\bigr)(0^+)=m_\alpha,
\]
and the solution is
\begin{equation}\label{eq:methane-rl-sol}
m(t)=m_\alpha t^{\alpha-1}E_{\alpha,\alpha}(-\lambda_\alpha t^\alpha)
+E t^\alpha E_{\alpha,\alpha+1}(-\lambda_\alpha t^\alpha).
\end{equation}
The two formulations are therefore not interchangeable at the initial time. Since the atmospheric burden is ordinarily prescribed by a finite classical value $m(0)=m_0$, the regularized/Caputo form \eqref{eq:methane-caputo} is the formulation consistent with the transform calculation and solution actually used in \cite{chishtie2026}. This distinction should be kept explicit when the model is used for inverse identification.

The stochastic motivation remains significant: the Monte Carlo construction in \cite{chishtie2026} reproduces the Mittag--Leffler survival law associated with subordination, and the authors report small mean numerical errors over the tested range of orders. This validates the numerical realization of the assumed CTRW mechanism; it does not by itself establish that global atmospheric methane observations require a unique fractional order. Seasonality, spatial exchange among reservoirs, uncertainty in OH chemistry, and temporal aggregation can all generate nonexponential effective responses.

A natural extension is distributed-order dynamics. For example,
\begin{equation}\label{eq:distributed-order}
\int_0^1\mu(\alpha)\,\CD{\alpha}m(t)\,d\alpha=F(m,t),
\end{equation}
where a physically interpretable weight should satisfy at least $\mu\in L^1(0,1)$, $\mu(\alpha)\ge0$ almost everywhere, and $\mu\not\equiv0$; one may additionally normalize $\int_0^1\mu(\alpha)\,d\alpha=1$ when $\mu$ is interpreted as a probability density of orders. 

\section{An electrochemical case study: audit before fractional interpretation}\label{sec:case}

A deviation from a simple classical formula is not, by itself, evidence for fractional dynamics. This section re-examines published electrochemical data from the laboratory of one of the authors \cite{rajabova2026}. The objective is not to reclassify the experiment as fractional, but to show why unit consistency, the choice of current law, and parameter identifiability must be checked before an anomalous exponent is given a physical interpretation.

\textit{1. Experimental setting and the observation model.} Paracetamol oxidation was studied on a glassy-carbon electrode of diameter 2~mm, modified with plant-derived gold nanoparticles and polyvinylpyrrolidone (AuNPs/PVP/GCE) \cite{rajabova2026}. The experiment combines several processes,
\[
\text{bulk transport}+\text{adsorption/interfacial structure}+\text{electron transfer}+\text{follow-up chemistry}.
\]
The published work describes paracetamol oxidation as an irreversible two-electron, two-proton process producing $N$-acetyl-$p$-benzoquinone imine, with diffusion and surface effects both contributing to the voltammetric response.

For a reversible diffusion-controlled couple at approximately 25\,$^\circ$C, the Randles--Sevcik relation is
\begin{equation}\label{eq:rs}
I_p=2.69\times10^5\,n^{3/2}ACD^{1/2}v^{1/2},
\end{equation}
when $I_p$ is in A, $A$ in cm$^2$, $C$ in mol\,cm$^{-3}$, $D$ in cm$^2$\,s$^{-1}$, and $v$ in V\,s$^{-1}$. The numerical coefficient is unit- and temperature-dependent. For an irreversible charge-transfer process the peak-current relation contains the transfer coefficient and the number of electrons involved in the rate-determining step; a commonly used 25\,$^\circ$C form is 
\begin{equation}\label{eq:rs-irr}
I_p=2.99\times10^5\,n(\alpha_a n_a)^{1/2}ACD^{1/2}v^{1/2}.
\end{equation}
Thus, even before any fractional model is considered, the electrochemical mechanism determines which observation formula is admissible.

In \cite{rajabova2026}, the reversible form \eqref{eq:rs} was first used with the [Fe(CN)$_6$]$^{3-/4-}$ probe to estimate an electroactive area $A=0.110$~cm$^2$, using $D=7.6\times10^{-6}$~cm$^2$\,s$^{-1}$. The same area was then inserted into a pointwise calculation of a paracetamol diffusion coefficient. Because the paracetamol process is described in the same paper as irreversible and adsorption-influenced, this second use of \eqref{eq:rs} is already model-dependent.

\textit{2. Arithmetic and unit audit of the reported diffusion coefficient.} A more basic difficulty appears when the published numbers are substituted back into the stated formula. Table~3 of \cite{rajabova2026} reports, for example,
\[
v=0.02\ {\rm V\,s^{-1}},\quad n=2,\quad A=0.110\ {\rm cm^2},\quad
C=0.05\ \mu{\rm M}=5.0\times10^{-11}\ {\rm mol\,cm^{-3}},
\]
and lists $D=4.01\times10^{-3}$~cm$^2$\,s$^{-1}$. The same row lists the numerical current as $3.03\times10^{-5}$ under a column headed ``$I_{pa}$ ($\mu$A)'', whereas Fig.~14 of that paper labels the current axis in amperes and shows currents of order $10^{-5}$~A. The two unit conventions are therefore inconsistent.

The discrepancy is not small. Substitution of the listed $D=4.01\times10^{-3}$~cm$^2$\,s$^{-1}$ into \eqref{eq:rs}, with the other reported quantities, predicts
\begin{equation}\label{eq:current-audit}
I_p\approx3.75\times10^{-8}\ {\rm A}=3.75\times10^{-2}\ \mu{\rm A},
\end{equation}
about $8.1\times10^2$ times smaller than the $3.03\times10^{-5}$~A scale visible in Fig.~14. If the table heading ``$\mu$A'' is taken literally, the listed number corresponds instead to $3.03\times10^{-11}$~A. Neither interpretation reproduces the reported $D$. Conversely, treating $3.03\times10^{-5}$ as amperes in \eqref{eq:rs} gives a diffusion coefficient of order $10^3$~cm$^2$\,s$^{-1}$, whereas treating it as microamperes gives a value of order $10^{-9}$~cm$^2$\,s$^{-1}$. These orders of magnitude are mutually incompatible.

The published value $D_{\rm app}=3.97\times10^{-3}$~cm$^2$\,s$^{-1}$ is therefore retained here only as a \emph{reported apparent parameter}; it is not used as a quantitatively validated molecular diffusion coefficient. A Wilke--Chang estimate of order $10^{-5}$~cm$^2$\,s$^{-1}$ provides a useful reference scale for a small organic solute in water \cite{wilke1955}, but no numerical ratio between that estimate and the reported electrochemical $D_{\rm app}$ is interpreted physically until the current units and the irreversible-current model are reconciled. This audit strengthens, rather than weakens, the methodological point:
\begin{center}
\emph{failure of a classical parameter calculation must first trigger a check of units, assumptions, and identifiability; it does not imply fractional transport.}
\end{center}

Table~\ref{tab:lab} records the quantities that can be checked directly from the publication and distinguishes them from re-analysed graphical quantities.

\begin{table}[ht]
\centering\small
\caption{Audit of transport-related quantities reported for AuNPs/PVP/GCE in \cite{rajabova2026}. ``Not reported'' means that an uncertainty or raw-data quantity needed for propagation of uncertainty was not supplied in the source.}
\label{tab:lab}
\begin{tabular}{>{\raggedright\arraybackslash}p{0.34\textwidth}>{\raggedright\arraybackslash}p{0.27\textwidth}>{\raggedright\arraybackslash}p{0.31\textwidth}}
\toprule
\textbf{Quantity} & \textbf{Reported value} & \textbf{Audit note}\\
\midrule
Probe diffusion coefficient & $7.6\times10^{-6}$~cm$^2$\,s$^{-1}$ & Assumed in area calibration; uncertainty not reported\\
Electroactive area & $0.110$~cm$^2$ & Obtained from the probe using \eqref{eq:rs}; transferred to paracetamol calculation\\
Paracetamol concentration & $0.05~\mu$M $=5.0\times10^{-11}$~mol\,cm$^{-3}$ & Unit conversion is unambiguous\\
Paracetamol scan rates & $0.02$--$0.06$~V\,s$^{-1}$ & Five points; highest-rate current is nonmonotone\\
Paracetamol current, first row & numerical value $3.03\times10^{-5}$ & Table header says $\mu$A; published voltammogram axis says A\\
Reported paracetamol $D$ & $(2.75$--$4.51)\times10^{-3}$~cm$^2$\,s$^{-1}$ & Cannot be reproduced from the stated current units and \eqref{eq:rs}\\
Temperature & 293.15~K in Table~3 & Coefficient $2.69\times10^5$ is the standard 25\,$^\circ$C form; the small temperature correction does not resolve the unit discrepancy\\
Uncertainty in $A$, $I_p$, $D$ & Not reported & Prevents a defensible uncertainty propagation for $D$\\
\bottomrule
\end{tabular}
\end{table}

\begin{figure}[t]
\centering
\includegraphics[width=\textwidth]{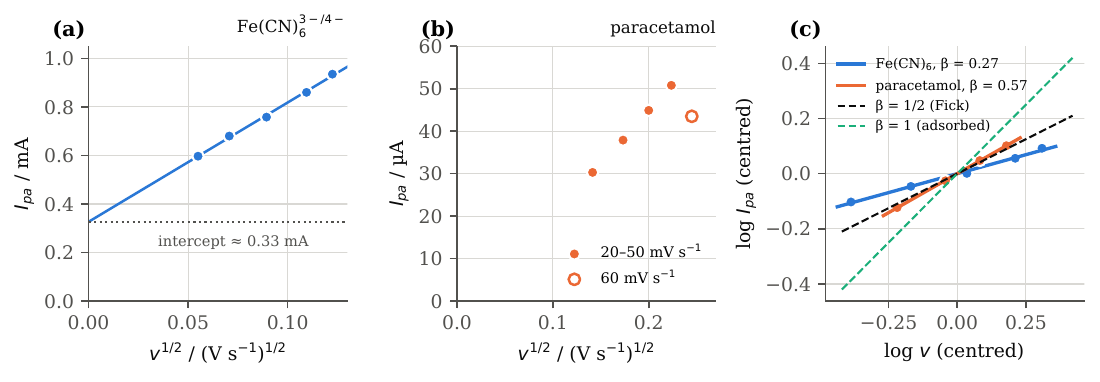}
\caption{Re-analysis of scan-rate data from \cite{rajabova2026}. (a) Hexacyanoferrate probe on AuNPs/PVP/GCE, digitized from the published figure, with an unconstrained regression in $v^{1/2}$; the nonzero intercept is shown explicitly. (b) Paracetamol currents on the ampere scale of the published voltammogram; the open symbol at 60~mV\,s$^{-1}$ is excluded from the four-point power-law fit. (c) Centered log--log slopes. Because a multiplicative current-unit conversion cancels in a logarithmic slope, the exponent is not affected by the A/$\mu$A ambiguity, although its physical interpretation remains model-dependent.}
\label{fig:scanrate}
\end{figure}

\textit{3. Scan-rate scaling: useful, but conditional.} The logarithmic slope
\[
\beta=\frac{d\log I_p}{d\log v}
\]
is unit-invariant and is therefore less sensitive to the current-unit inconsistency than an absolute diffusion-coefficient calculation. For ideal semi-infinite Fickian diffusion, $\beta=1/2$; for an ideal surface-confined redox species, $\beta=1$. Figure~\ref{fig:scanrate} re-plots the published scan-rate series. The paracetamol current values in panel (b) are read on the ampere scale of the published voltammogram; multiplying all currents by a constant unit conversion would not change the fitted logarithmic exponent.

\emph{Hexacyanoferrate probe.} Digitization of the published probe plot gives an unconstrained fit with a substantial positive intercept, and a centered log--log slope of approximately $\beta\approx0.27$. A high $R^2$ for a straight line in $I_p$ versus $v^{1/2}$ does not by itself enforce the zero-intercept scaling required by an ideal Randles--Sevcik law. Background current, capacitive current, or other surface contributions can change an apparent exponent over a short scan-rate window.

\emph{Paracetamol.} Using the first four points (20--50~mV\,s$^{-1}$) gives $\beta\approx0.57$; the 60~mV\,s$^{-1}$ point is lower, consistent with the source's discussion of peak broadening, capacitive contributions, and adsorption of oxidation products. The scan-rate interval spans only a factor of $2.5$ for the four-point fit, so the exponent should be regarded as descriptive rather than as a precisely identified asymptotic power law.

The two flux closures of Section~4.2 show why these slopes cannot be converted automatically into fractional orders. The identities $\beta=\alpha/2$ and $\beta=1-\alpha/2$ hold only under their respective ideal planar diffusion models and current definitions. The present paracetamol process is irreversible and adsorption-influenced, so neither equality is an estimator of $\alpha$ without an electrode-kinetic model. The slopes are therefore used here only to demonstrate that the data do not determine a unique transport mechanism.

\textit{4. Impedance data.} The published impedance spectra were recorded with 5~mM [Fe(CN)$_6$]$^{3-/4-}$/0.1~M KCl and fitted with the Randles circuit $R_s(C_{dl}(R_{ct}W))$ \cite{rajabova2026}. The reported values are $R_s\approx314$~$\Omega$ and $R_{ct}=72.4$ and $99.3$~k$\Omega$ for the bare and modified electrodes. Figure~\ref{fig:eis}(b) digitizes the published Nyquist curves and normalizes each by its reported charge-transfer resistance. The last few graphical points suggest a low-frequency line steeper than the $45^\circ$ ideal Warburg direction, with an angle of roughly $63^\circ$. Because raw impedance values and uncertainties were not available for the present re-analysis, this angle is an approximate graphical diagnostic rather than a fitted material parameter.

A non-Warburg slope is not mechanism-specific. It can arise from finite-length diffusion, rough or porous geometry, distributed interfacial relaxation, or anomalous bulk transport \cite{bisquert2001}. Moreover, Section~4.3 shows that even within an anomalous-diffusion interpretation the mapping from impedance angle to transport order depends on the constitutive current closure. The impedance spectrum therefore motivates model discrimination; it does not identify a fractional exponent by itself.

\textit{5. Structural identifiability survives the unit audit.} The unit inconsistency prevents a reliable numerical use of the published $D_{\rm app}$, but it does not affect the structural observation that a Randles--Sevcik-type measurement identifies a product of parameters. If
\begin{equation}\label{eq:ident-case}
Y=ACD^{1/2},
\end{equation}
then for every $c>0$ the pairs $(A,D)$ and $(cA,D/c^2)$ give the same $Y$ when $C$ is fixed. No improvement in numerical optimization can separate $A$ and $D$ without an independent observation. Figure~\ref{fig:ident} illustrates this geometry in normalized variables, without using the disputed numerical diffusion coefficient.

\textit{6. What a convincing test would require.} A discriminating experimental program should compare several model classes under independent perturbations: a classical Fick/Butler--Volmer model with the appropriate reversible or irreversible peak-current relation; an interfacial model with adsorption, distributed capacitance, porous geometry, or finite-length diffusion; and one or more fractional transport models with explicitly stated flux closures. Useful data would include cyclic voltammetry over a substantially wider scan-rate range with blank-current subtraction and replicate uncertainties; chronoamperometry, for which the two ideal fractional closures predict different temporal exponents; impedance spectra at several DC potentials and film loadings with Kramers--Kronig validation; viscosity variation to perturb bulk transport; and independent area calibration with propagated uncertainty. A physically meaningful fractional order should be stable across independent observables and experimental perturbations, not merely improve one in-sample fit.

The case study therefore supports a conservative inference: the published data contain model-selection questions worth studying, but they do not presently establish fractional transport. The first mathematical task is to make the classical observation operator, units, and identifiable parameter combinations internally consistent.

\section{Well-posedness before numerical computation}

Fractional chemical models are frequently explored through numerical trajectories. Such trajectories can be scientifically interpreted only after the underlying continuous problem is specified well enough to determine what the computation is approximating. Consider abstractly
\begin{equation}\label{eq:abstract}
\CD{\alpha}u(t)=F(t,u(t)),\qquad u(0)=u_0 .
\end{equation}
At least five logically distinct questions arise.
\begin{enumerate}[label=(\roman*)]
\item \textit{Existence.} Does a solution satisfying the initial and boundary conditions exist in the intended class?
\item \textit{Uniqueness.} Is the solution determined uniquely by the data?
\item \textit{Physical admissibility.} Are concentrations nonnegative and are required mass, charge, or stoichiometric constraints preserved?
\item \textit{Continuous dependence.} Do small perturbations of data and parameters produce controlled changes in the solution?
\item \textit{Convergence of the numerical scheme.} Does the discrete approximation approach the continuous solution as the discretization is refined?
\end{enumerate}

For \eqref{eq:abstract} with $F$ continuous and Lipschitz in $u$, existence and uniqueness of a continuous solution follow from the equivalent Volterra integral equation, and multi-order systems can be treated by fixed-point arguments; stability notions for fractional systems differ in an essential way from the integer-order case. For chemical models the difficulty usually lies in the nonlinear, coupled or boundary-value setting. Even in the linear case the solution is typically not smooth at $t=0$: $u'(t)\sim t^{\alpha-1}$, so standard schemes on uniform meshes lose their nominal order of accuracy, which motivates graded meshes.

If a nonlinear fractional equation has several solutions, a numerical algorithm may converge to one branch because of the chosen discretization, initialization, or regularization. A smooth graph does not establish that this branch is unique or physically distinguished. Conversely, uniqueness of the continuous problem does not imply that an arbitrary discretization converges to it. Therefore
\begin{center}
\textbf{numerical solvability does not replace mathematical well-posedness.}
\end{center}

Recent chemical applications exhibit a wide range of mathematical completeness. Some enzyme and pyrolysis papers establish existence, uniqueness, boundedness, or Ulam-type stability before computation \cite{agilan2025,zulqarnain2026,pathak2026}; several reaction--diffusion papers provide explicit stability and convergence orders for discretizations \cite{roul2024,chawla2025}; other works are primarily phenomenological or computational. These are different levels of analysis and should not be conflated.

\section{Mathematical and physical status of representative model classes}

Table~\ref{tab:status} summarizes representative classes discussed above. ``Explicit representation'' means a closed formula or transform representation that determines the solution uniquely in the natural class. The final column identifies the principal mathematical or experimental issue that remains after the cited results are taken into account.

{\small
\setlength{\LTcapwidth}{\textwidth}
\begin{longtable}{>{\raggedright\arraybackslash}p{0.135\textwidth}>{\centering\arraybackslash}p{0.075\textwidth}>{\raggedright\arraybackslash}p{0.19\textwidth}>{\raggedright\arraybackslash}p{0.25\textwidth}>{\raggedright\arraybackslash}p{0.205\textwidth}}
\caption{Mathematical and physical status of representative fractional models in chemistry. Evidence codes: D---the nonlocal law is derived from a transport, stochastic, constitutive, or reduction mechanism; E---the model is confronted with measured chemical or physicochemical data; M---the principal contribution is mathematical or computational analysis of a chemically motivated model. Combinations indicate that more than one type of evidence is present.}\label{tab:status}\\
\toprule
\textbf{Model class} & \textbf{Basis} & \textbf{Origin of nonlocality} & \textbf{Established structure} & \textbf{Main remaining issue}\\
\midrule\endfirsthead
\toprule
\textbf{Model class} & \textbf{Basis} & \textbf{Origin of nonlocality} & \textbf{Established structure} & \textbf{Main remaining issue}\\
\midrule\endhead
\bottomrule\endfoot
Semi-integral electroanalysis & D/E & Elimination of space from classical semi-infinite diffusion & Explicit half-order Dirichlet-to-Neumann relation; classical well-posedness is inherited; semi-integral electroanalysis was subsequently tested experimentally & Extension to finite, curved, heterogeneous, and nonlinear interfaces\\\addlinespace[3pt]
Fractional voltammetry & D/M & Time-fractional transport with a specified current closure & Explicit transform/Wright representations; current scaling follows from the chosen closure and ideal electrochemical boundary conditions & Experimental discrimination between gradient and conservative memory-flux closures\\\addlinespace[3pt]
Anomalous impedance and PNP & D/E/M & Bulk anomalous transport, distributed order, or fractional boundary response & Generalized Warburg/CPE-type frequency laws; several formulations have been fitted to impedance data & Separation of bulk and interfacial memory; closure dependence; passivity and identifiability\\\addlinespace[3pt]
Fractional pseudo-first-order adsorption & E/M & Caputo relaxation used as an effective memory law & Explicit Mittag--Leffler solution; positivity and monotone approach to equilibrium in the linear model & Mechanistic identification versus heterogeneous classical alternatives\\\addlinespace[3pt]
SPR adsorption & E & Fractional one-to-one binding fitted to measured sensorgrams & Reported fitting error more than one order of magnitude below the integer one-to-one model & Comparison with mass-transfer, multisite, and multistate classical models; out-of-sample prediction\\\addlinespace[3pt]
Barrier--diffusion adsorption & D/M & Anomalous bulk flux followed by elimination of spatial variables & Conservation built into the formulation; linear Henry case has explicit representations & Global nonlinear well-posedness, positivity, and convergence for nonlinear isotherms\\\addlinespace[3pt]
Methanol in porous solids & E/D & Measured non-Fickian transport associated with trapping/adsorption in porous media & Fractional transport laws compared across related porous materials; fitted exponent changes with material structure & Separation of geometry, adsorption energetics, and the stochastic waiting-time mechanism\\\addlinespace[3pt]
Reaction--subdiffusion & D/M & CTRW competition between waiting and reaction clocks & Microscopic derivations show that memory generally enters the effective reaction term as well as transport & General nonlinear networks, complex boundaries, positivity, and stoichiometric invariants\\\addlinespace[3pt]
Conservation-consistent fractional kinetics & M & Fractional reaction-rate structure constrained by stoichiometry & Mass-balance-preserving constructions for nonlinear reactions & General network theory and experimental validation\\\addlinespace[3pt]
Oscillatory reactions & M & Mostly phenomenological Caputo generalizations; one cited study analyzes a discrete map derived from such a system & Continuous-time studies provide stability/numerical analyses; period-doubling and Neimark--Sacker results of Berkal--Almatrafi concern the discrete-time map & Experimental identification of the order and rigorous separation of continuous and discrete dynamics\\\addlinespace[3pt]
Enzyme models & M & Constant/variable memory, delays, cooperative kinetics & In several proposed models existence, uniqueness, boundedness, or stability is established & Mechanistic link to conformational statistics; structural identifiability; validation against independent biochemical data\\\addlinespace[3pt]
Pyrolysis & E/M & Phenomenological thermal memory in non-isothermal kinetics & Fractional models fitted to thermal-decomposition data; some analyses include positivity, boundedness, or existence/uniqueness & Cross-validation across heating protocols and comparison with distributed activation-energy models\\\addlinespace[3pt]
Polymers/CPE & D/E & Broad distributions of molecular/interfacial relaxation times & Fractional constitutive laws reproduce broadband relaxation; positive spectral representations are available for important kernels & Uniqueness of relaxation spectrum from finite-band data; passivity and parameter identifiability\\\addlinespace[3pt]
Atmospheric methane & D/M & CTRW/subordination and heavy-tailed residence times & Explicit linear Mittag--Leffler relaxation and Monte Carlo realization of the assumed stochastic mechanism; the cited paper's RL label must be distinguished from the regularized transform actually used & Stable recovery of an order from atmospheric observations under forcing and aggregation\\\addlinespace[3pt]
Fractional quantum chemistry & M & Nonlocal quantum operators proposed at the model level & Fractional HF/KS-DFT and related operator frameworks have been formulated & Operator domains, variational stability, convergence, and calibration against standard electronic-structure data\\\addlinespace[3pt]
Authors' laboratory data \cite{rajabova2026} & E & Classical electrochemical observation model with diffusion, adsorption, and interfacial effects & The published paracetamol diffusion coefficient cannot be reproduced from the displayed current units and stated Randles--Sevcik formula; scan-rate and impedance plots contain nonclassical features but do not identify a unique mechanism & Reconcile units and irreversible kinetics; then compare classical interfacial and fractional models with propagated uncertainty and identifiability analysis\\\addlinespace[3pt]
\end{longtable}
}

The table separates two statements that are often confused: \emph{the mathematical problem has a unique solution} and \emph{the chemical system requires a fractional model}. The first concerns the internal validity of a model; the second requires derivation, discriminatory experiments, or both.

\section{Nonnegativity and conservation laws as chemical constraints}

Existence in an abstract function space is not sufficient for a concentration model. If $C(x,t)$ denotes concentration, a negative value usually has no physical meaning. Preservation of the positive cone is therefore part of physical admissibility.

Likewise, chemical reactions inherit stoichiometric invariants. If the reaction network conserves a chemical component, then for suitable coefficients $a_j$, $\sum_ja_jC_j(t)=\mathrm{const}$. In matrix form, if $\dot c=S\,r(c)$ and $\ell^TS=0$, then $\ell^Tc(t)$ is conserved. A fractional generalization should preserve the same invariant unless the model explicitly adds a source or sink.

This requirement is especially delicate when different fractional operators are assigned independently to different species. Even a mathematically solvable system can then violate the stoichiometric structure of the chemistry. Conservation-consistent formulations, such as the rate-based construction in \cite{barajas2025}, avoid this problem by retaining a common reaction flux across component balances. 

Transport models provide a parallel lesson. The safest architecture is to preserve the local balance $\partial_tC=-\nabla\cdot J$ and place anomalous behavior in a constitutive relation for $J$. This is the structure of the barrier--diffusion adsorption model \eqref{eq:bb1}. It separates the universal balance law from the material-specific transport law and makes mass conservation transparent. It differs from the gradient-based current closure used for the Mirceski--Tomovski model; the two formulations coincide in their classical limit $\alpha=1$, but for $\alpha<1$ they predict different current exponents, hence the constitutive flux law must always be stated explicitly.

\section{Dimensional consistency}

A fractional derivative has a physical dimension. If $D_t^\alpha$ is of order $\alpha$, then $[D_t^\alpha]=T^{-\alpha}$. Thus $\CD{\alpha}u=-\lambda_\alpha u$ requires $[\lambda_\alpha]=T^{-\alpha}$, not $T^{-1}$. Likewise, in a time-fractional diffusion equation the generalized diffusion coefficient typically has dimension $L^2T^{-\alpha}$, so it cannot be compared numerically with an ordinary $D$ unless a time scale is specified. Under the scaling $t=T_0\tilde t$, $D_t^\alpha=T_0^{-\alpha}D_{\tilde t}^\alpha$.

Nondimensionalization is therefore a useful diagnostic rather than a cosmetic change of variables. The anomalous adsorption model explicitly carries the fractional dimensions through the reduction before passing to dimensionless variables \cite{bazhlekov2026}. Similar care is taken in the batch-reactor formulation of Qureshi and Aziz \cite{qureshi2020}. Dimensional consistency cannot validate a physical model, but dimensional inconsistency is strong evidence that a fractional operator has been inserted without sufficient care. In variable-order equations the issue is even more delicate because a coefficient that compensates $T^{-\alpha(t)}$ must be defined consistently as the order changes.

\section{Parameter identifiability and experimental design}

Well-posedness of the forward problem does not imply that its parameters can be recovered uniquely from an experiment. Suppose the measured quantity has the form
\begin{equation}\label{eq:ident}
Y=AD^{1/2}.
\end{equation}
Even exact knowledge of $Y$ cannot determine $A$ and $D$ separately: for every $c>0$, the pairs $(A,D)$ and $(cA,D/c^2)$ produce the same observation. This is structural nonidentifiability; no improvement in measurement precision or numerical optimization can remove it without additional information.

\begin{figure}[t]
\centering
\includegraphics[width=0.62\textwidth]{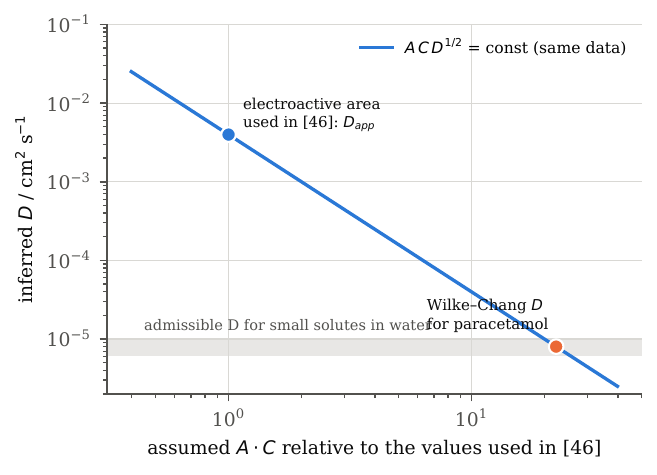}
\caption{Structural non-identifiability of an observation $Y=P D^{1/2}$ with $P=AC$. In normalized variables all points on $D/D_0=(P/P_0)^{-2}$ give exactly the same value of $Y$. The diagram is deliberately independent of the disputed absolute current units in the paracetamol data of \cite{rajabova2026}; it illustrates the observation geometry rather than a fitted diffusion coefficient.}
\label{fig:ident}
\end{figure}

The same geometry applies to a Randles--Sevcik-type observation, for which the measured current is proportional to $ACD^{1/2}$ when the remaining factors are fixed. The unit inconsistency identified in Section~\ref{sec:case} prevents a reliable numerical placement of the paracetamol experiment on such a curve, but it does not alter the structural statement. Figure~\ref{fig:ident} therefore shows the relation in normalized variables. If $P=AC$ and $Y=PD^{1/2}$ is fixed, then
\[
\frac{D}{D_0}=\left(\frac{P}{P_0}\right)^{-2}.
\]
Every point on this curve produces the same observation $Y$; an independent measurement of $P$ or a second observable with different parameter dependence is required to determine $D$.

Fractional models introduce analogous correlations. The order $\alpha$ may trade off against a generalized diffusion coefficient, a reaction constant, a CPE amplitude, an interfacial capacitance, or a delay. Heavy-metal sorption provides a useful example in which joint confidence regions are used to test whether the fitted order is statistically separated from the integer value \cite{friesen2015}. Similar analyses should become standard in chemical applications of fractional calculus.

Experimental design should therefore combine observables with different parameter dependence. In electrochemistry, cyclic voltammetry, chronoamperometry, and impedance spectroscopy probe different combinations of transport and interface parameters. In adsorption, multiple concentrations and temperatures can distinguish intrinsic binding from mass transfer. In enzyme kinetics, substrate variation and independent measurements of conformational times can help separate kinetic constants from memory parameters. In atmospheric chemistry, several temporal resolutions can test whether the inferred order is an artifact of aggregation.

One must also distinguish parameter identifiability from model identifiability. Even if $\alpha$ is uniquely estimated inside one fractional model, a heterogeneous classical model may give nearly the same output. The inverse problem may therefore require a discrete choice among model classes as well as continuous parameter estimation within each class.

\section{What constitutes evidence for a fractional chemical model?}

The literature suggests a hierarchy of evidential strength.
\begin{enumerate}
\item \textit{Derivation from a more fundamental model.} The half-order boundary operator generated by semi-infinite classical diffusion is the clearest example.
\item \textit{Microscopic stochastic limit.} A CTRW with heavy-tailed waiting times can lead to a fractional transport or reaction--transport equation \cite{metzler2000}.
\item \textit{Independent scaling law.} Mean-square displacement, current decay, impedance slope, or first-passage behavior should follow a power law predicted before fitting. In electrochemistry the relation between the observed exponent and $\alpha$ is closure-dependent.
\item \textit{Transferable prediction across several experiments.} The same order should predict data not used in calibration, preferably under changes in concentration, scan rate, temperature, or geometry.
\item \textit{Statistical identifiability.} Confidence regions should exclude the classical value when appropriate, and the order should remain stable when the observation window is varied.
\item \textit{Improved fit to a single data set.} This is useful but weaker, because mixtures of classical exponentials, heterogeneous sites, or unmodeled transport can produce similar curves.
\item \textit{Formal replacement $d/dt\mapsto D_t^\alpha$.} Without a mechanism, a discriminatory observable, or a new verified prediction, this is the weakest justification.
\end{enumerate}
This hierarchy does not restrict the use of fractional calculus. It clarifies which conclusions can legitimately be drawn from each type of evidence and helps separate mathematical exploration from physical identification.

\section{Open problems and directions for fractional analysis in chemistry}

The preceding models lead to a set of concrete mathematical problems in which advances in fractional analysis could have a direct chemical impact.

\textit{Open Problem 1: global well-posedness of nonlinear barrier--diffusion adsorption.}
Consider the reduced two-sided adsorption system
\begin{equation}\label{eq:op1}
\frac{d\sigma^\pm}{dt}+a_\pm\RL{\gamma_\pm/2}\sigma^\pm
=b_\pm-c_\pm f(\sigma^++\sigma^-),\qquad \sigma^\pm(0)=0,
\end{equation}
where $0<\gamma_\pm<2$, $a_\pm,b_\pm,c_\pm>0$, and $f$ is the nonlinear isotherm term. A natural target is a theorem in the class $\sigma^\pm\in AC_{\rm loc}([0,\infty))$ with $\RL{\gamma_\pm/2}\sigma^\pm\in L^1_{\rm loc}$, under assumptions such as $f(0)=0$, local Lipschitz continuity and monotonicity on the physically admissible interval $[0,\sigma_\infty)$. Prove global existence and uniqueness, continuous dependence, and invariance of $0\le \sigma^+(t)+\sigma^-(t)<\sigma_\infty$; determine conditions for monotone versus nonmonotone adsorption. The Henry case $f(\sigma)=k\sigma$ is linear and already has explicit generalized Mittag--Leffler representations in \cite{bazhlekov2026}; the global nonlinear result for Volmer and van der Waals isotherms is the open part.

\textit{Open Problem 2: rigorous convergence of nonlinear memory solvers.}
Integrating \eqref{eq:op1} and using $\sigma^\pm(0)=0$ gives the Volterra form
\begin{equation}\label{eq:op2}
\sigma^\pm(t)=b_\pm t-c_\pm\int_0^t f(\sigma^+(s)+\sigma^-(s))\,ds
-a_\pm\Ii{1-\gamma_\pm/2}\sigma^\pm(t).
\end{equation}
For this equation, establish an error theory for the predictor--corrector discretization used in the nonlinear adsorption literature and for fast-convolution variants. The theorem should state the regularity assumed near $t=0$, give an explicit convergence order on uniform or graded meshes, and prove that the fully discrete solution remains in the physical invariant interval whenever the continuous solution does. A further target is a fast implementation with near-linear or $O(N\log N)$ work for $N$ time levels without destroying positivity or the proven error bound. Existing computations in \cite{bazhlekov2026} demonstrate the method numerically but do not provide this complete convergence-and-invariance theory.

\textit{Open Problem 3: inverse recovery of order, flux law, and transport parameters.} Investigate uniqueness and stability for recovery of $(\alpha,d_\alpha,k_0,A,\dots)$ from combined electrochemical observations. In particular, determine when data can distinguish the gradient closure from a conservative memory-flux closure rather than merely estimate $\alpha$ inside one assumed model.Without an identifiability analysis the learned order need not be unique. The electrochemical setting differs in that the observation is a boundary flux and the parameter vector also contains the area and kinetic constants. This problem is directly motivated by the scan-rate and impedance data of Section~\ref{sec:case}.

\textit{Open Problem 4: distributed-order and variable-order laws derived from mechanisms.} Derive distributed-order or variable-order equations from explicit mixtures of residence times, evolving pore structures, conformational states, or environmental regimes. Develop inverse theories for the order distribution or order function and identify conditions under which they are distinguishable from slowly varying classical parameters.

\textit{Open Problem 5: classical heterogeneous versus genuinely fractional models.} Construct statistically and mathematically controlled model-selection procedures that compare fractional equations with mixtures of exponentials, distributed RC networks, multisite adsorption, multicompartment kinetics, and heterogeneous classical transport. Preference for a fractional model should be based on predictive performance on held-out data, not only on in-sample residual reduction.

\textit{Open Problem 6: a direct--inverse--numerical--experimental theory.} For a small number of chemically well-motivated systems, complete the chain
\begin{align*}
&\text{well-posed forward problem}\to\text{stable inverse problem}\\
&\quad\to\text{convergent numerical method}\to\text{experimental validation}.
\end{align*}
At present these four components are usually studied separately.

\textit{Open Problem 7: physically consistent nonlinear reaction--subdiffusion networks.}
Let $S\in\mathbb R^{m\times r}$ be a stoichiometric matrix and let $r(c)$ be a locally Lipschitz, quasipositive mass-action rate vector on $\mathbb R_+^m$. A satisfactory reaction--subdiffusion model should be derived from a specified CTRW or semi-Markov mechanism, reduce to the known first-order reaction--subdiffusion equations of \cite{seki2003a,seki2003b} in the corresponding linear limit, and preserve every stoichiometric invariant $\ell^Tc$ with $\ell^TS=0$. The open task is to obtain a closed nonlinear evolution law for general $S r(c)$ when reactions may occur during trapping, and then prove positivity, conservation of all such invariants, and global well-posedness for nonnegative initial data. In particular, the desired law should make precise when the frequently written additive model
\[
\partial_t c_i=d_i\RL{1-\alpha_i}\Delta c_i+(Sr(c))_i
\]
is valid and when the reaction term itself must carry the memory generated by the waiting-time process.

\textit{Open Problem 8: fractional transport in porous catalysts and reactor design.} Relate microscopic trapping and adsorption statistics to effective fractional parameters and then to reactor-level quantities such as the Thiele modulus, effectiveness factor, optimal particle size, and conversion. Establish whether parameters inferred from transient uptake predict steady and periodically forced reactor behavior.

\textit{Open Problem 9: passivity and spectral representation of impedance and polymer models.} Characterize those fractional impedance and dielectric kernels that correspond to causal passive systems. For the CPE and the Cole--Cole law the equivalence with positive distributions of relaxation times is classical: their relaxation functions are completely monotone and admit explicit positive spectra \cite{gorenflo2020}. Determine when a more general fractional constitutive law has a unique positive distribution of ordinary relaxation times and quantify the nonuniqueness caused by a finite frequency window.

\textit{Open Problem 10: experimental identification in oscillatory and enzyme kinetics.} Design experiments in which the same fractional order must predict several independent observables. For Brusselator and BZ systems these may include periods, amplitudes, phase response, and wave speed; for enzyme kinetics they may include multiple substrate concentrations, residence-time distributions, and single-molecule conformational statistics.

\textit{Open Problem 11: conservation-preserving fast numerical methods.}
As a concrete prototype, consider a spatially discretized fractional reaction--diffusion system
\begin{equation}\label{eq:op11}
\CD{\alpha}U(t)+AU(t)=F(U(t)),\qquad U(0)=U_0\ge0,
\end{equation}
where $0<\alpha<1$, $A$ is an $M$-matrix (or, before discretization, generates a positive diffusion semigroup), $F$ is locally Lipschitz and quasipositive, and a prescribed vector $\ell$ satisfies $\ell^TA=0$ and $\ell^TF(U)=0$ so that $\ell^TU$ is an invariant. Develop a high-order fast-memory scheme for \eqref{eq:op11} that (i) converges for the typical weak initial singularity $U'(t)=O(t^{\alpha-1})$, with an explicit rate on graded/adaptive meshes; (ii) preserves $U^n\ge0$ and the invariant $\ell^TU^n=\ell^TU^0$ exactly or with a quantified defect; and (iii) reduces the history cost from the direct quadratic convolution to near-linear work.

\section{Conclusion}

Fractional differential equations have a genuine role in chemical modeling when the nonlocal operator can be connected to a specific mechanism. Several mechanisms recur across apparently different applications. A half-order boundary operator emerges exactly when the spatial variable is eliminated from semi-infinite classical diffusion. Heavy-tailed waiting times generate time-fractional transport and, when reaction competes with trapping, modify the reaction term itself. Broad distributions of relaxation times produce CPE-like impedance, fractional Langmuir behavior, and non-Debye polymer relaxation. Anomalous bulk transport coupled to an interface generates reduced fractional adsorption equations while preserving local mass balance. In voltammetry and impedance, the relation between a measured exponent and the fractional order depends on the constitutive definition of the flux.

At the same time, the limitations of formal fractionalization remain important. In oscillatory reactions, enzyme models, and several environmental systems, the Caputo order is often introduced phenomenologically and then shown to alter numerical trajectories or stability. Such results may be mathematically interesting, but they do not establish that the underlying chemistry has a power-law memory. Physical identification requires either derivation from a microscopic or reduced model, or independent experiments that distinguish the fractional model from plausible classical heterogeneous alternatives.

Mathematical reliability is inseparable from that physical question. Existence, uniqueness, continuous dependence, nonnegativity, conservation laws, dimensional consistency, structural identifiability, and numerical convergence determine whether a calculated trajectory has a defensible interpretation. A model can be mathematically impeccable yet chemically inappropriate; conversely, an empirically successful fit can rest on an ill-posed or non-identifiable mathematical formulation. The strongest future work will connect both levels.

The electrochemical case study illustrates the point in practical form. A published diffusion-coefficient calculation is not reproducible from its stated current units and observation formula, while digitized scan-rate and impedance plots contain nonclassical slopes whose interpretation depends on background subtraction, adsorption, geometry, electrode kinetics, and the constitutive current law. The correct conclusion is therefore not that fractional transport has been detected, but that the observation model must be made internally consistent and then compared against competing classical and fractional mechanisms. Fractional analysis becomes scientifically useful when it helps design experiments that separate these alternatives rather than merely replacing one fitted exponent by another.

Chemistry offers a rich collection of multiscale systems with memory, trapping, interfaces, heterogeneous transport, and distributed relaxation. Fractional analysis provides tools for describing such phenomena, but its greatest contribution will come from models in which derivation, rigorous analysis, computation, inverse identification, and experiment are developed together.

\section*{Declaration of generative AI use}

During the preparation of this work the authors used ChatGPT (OpenAI),
large language model, to review the manuscript, to revise the wording of the text, to help write and to suggest and check
bibliographic references. After using this tool, the authors reviewed and edited the content, verified
the mathematical arguments, and take full responsibility for the content
of the publication. The tool is not listed as an author.

\bibliographystyle{amsplain}

\end{document}